\magnification=\magstep1
\input amstex
\documentstyle{amsppt}
\hoffset=.25truein
\hsize=6truein
\vsize=8.75truein

\topmatter
\centerline{\bf ON THE TRACE OF HECKE OPERATORS FOR}
\title
 MAASS FORMS FOR CONGRUENCE SUBGROUPS II
\endtitle
\abstract
  Let $E_\lambda$ be the Hilbert space spanned by the eigenfunctions of the
non-Euclidean Laplacian associated with a positive
discrete eigenvalue $\lambda$.
In this paper, the trace of Hecke operators $T_n$ acting on the
space $E_\lambda$ is computed for Hecke congruence subgroups $\Gamma_0(N)$
of non-square free level.  This extends the computation of Conrey-Li [2],
where only Hecke congruence subgroups $\Gamma_0(N)$ of square free
level $N$ were considered.
\endabstract
\author
 Xian-Jin Li
\endauthor
\address
Department of Mathematics, Brigham Young University, Provo, Utah 84602
\endaddress
\email
xianjin\@math.byu.edu
\endemail
\thanks
Research partially supported by the American Institute of Mathematics.
\endthanks
\endtopmatter
\document

\heading
1.   Introduction
\endheading

   Let $N>1$ be an integer, which is not square free.  Denote by
$\Gamma_0(N)$ the Hecke congruence subgroup of level $N$.
The non-Euclidean Laplacian  $\Delta$ on the upper half-plane
$\Cal H$ is given by
$$\Delta= -y^2\left(\frac{\partial ^2}{\partial x^2}+\frac{\partial
^2}{\partial y^2} \right).$$
Let $D$ be the fundamental domain of $\Gamma_0( N)$.   Eigenfunctions
of the discrete spectrum of $\Delta$ are nonzero
real-analytic solutions of the equation
$$\Delta\psi=\lambda\psi$$
such that $\psi(\gamma z)=\psi(z)$ for all $\gamma$ in $\Gamma_0( N)$ and
such that
$$\int_D |\psi(z)|^2dz<\infty$$
where $dz$ represents the Poincar\'e measure of the upper half-plane.

  Let $\frak a$ be a cusp of $\Gamma_0(N)$.  We denote its
 stabilizer by $\Gamma_{\frak a}$.
  An element $\sigma_{\frak a}$ in $PSL_2(\Bbb R)$
exists such that $\sigma_{\frak a}\infty=\frak a$ and
$\sigma_{\frak a}^{-1}\Gamma_{\frak a}\sigma_{\frak a}=\Gamma_\infty$.
If $\psi$ is an eigenfunction of the Laplacian associated
with a positive discrete eigenvalue
$\lambda$, then it has a Fourier expansion [6]
$$\psi(\sigma_{\frak a}z)=\sqrt y\sum_{m\neq 0}\rho_{\frak a}(m)
K_{i\kappa}(2\pi |m|y)e^{2m\pi ix}$$
at every cusp $\frak a$ of $\Gamma_0(N)$,
where $\kappa=\sqrt{\lambda-1/4}$ and $K_\nu(y)$ is given by
the formula $\S$6.32, [18]
$$ K_\nu(y)=\frac{2^\nu\Gamma(\nu+\frac 12)}{y^\nu\sqrt\pi}\int_0^\infty
\frac{\cos (yt)}{(1+t^2)^{\nu+\frac 12}}dt. \tag 1.1$$

   The Hecke operators $T_n$, $n=1, 2, \cdots$, $(n, N)=1$, which act
in the space of automorphic functions with respect to $\Gamma_0( N)$,
are defined by
$$\left(T_nf\right)(z)=\frac 1{\sqrt n}\sum_{ad=n, \,
0\leq b<d}f\left(\frac {az+b}d\right).$$
An orthonormal system of eigenfunctions of $\Delta$ exists [6] such
that each of them is an eigenfunction of all the Hecke operators.  We call
these eigenfunctions Maass-Hecke eigenfunctions.
Let $\lambda_j$, $j=1, 2, \cdots$, be an enumeration in increasing order
of all positive discrete eigenvalues of $\Delta$ for $\Gamma_0(N)$
with an eigenvalue of multiplicity $m$ appearing $m$ times, and let
$\kappa_j=\sqrt{\lambda_j-1/4}$.   If $\psi_j(z)$ is a Maass-Hecke
 eigenfunction of $\Delta$ associated with the $j$th eigenvalue
$\lambda_j $, then
$$\left(T_n\psi_j\right)(z)=\tau_j(n)\psi_j(z)$$
where $\rho_{j\infty}(m)=\rho_{j\infty}(d)\tau_j(n)$ if $m=dn$
with $n\geqslant 1$, $(n, dN)=1$.

  Let $E_\lambda$ be a Hilbert space of functions spanned by
the eigenfunctions of $\Delta$ associated with a positive
discrete eigenvalue
$\lambda$.  The Petersson inner product of the space is given by
$$\langle F(z), G(z)\rangle=\int_D F(z)\bar G(z)dz. \tag 1.2$$
The analogue for Maass forms of Eichler-Selberg's trace formula
[12], p.85 is obtained for congruence
subgroups $\Gamma_0(N)$ of square free level $N$ in Conrey-Li [2].
In this paper, the trace $\text{tr}_\lambda T_n$ of Hecke operators acting
on the space $E_\lambda$ is computed for congruence subgroups
$\Gamma_0(N)$ of non-square free $N$.  Some of this computation
is implicit in Hejhal [4][5].

    Denote by $h_d$ the class number of
indefinite rational quadratic forms with discriminant d.  Let
$$\epsilon_d=\frac{v_0+u_0\sqrt d}2, \tag 1.3$$
where the pair $(v_0, u_0)$ is the fundamental solution [10] of Pell's
equation $v^2-du^2=4$.  Denote by $\Omega$ the set of all the positive
integers $d$ such that $d\equiv 0$ or 1 (mod 4) and such that $d$ is not
a square of an integer.  Let $\mu(n)$ be the M\"obius function.

\proclaim{Theorem 1}    Let $N$ be a non-square free positive
integer, and let $n$ be a positive integer with $(n, N)=1$.
If
$$L_n(s)=\sum_{k|N}k^{1-2s}\sum_{m|{N\over k}}\sum_{n|k}{\mu(n)\over n}
\sum_{d\in \Omega}\sum_u \prod_{p^{2l}|(d, {N\over k})}p^l
\left({d\over mn}\right){h_d\ln \epsilon_d\over (du^2)^s} $$
for $\text{Re}\, s>1$, where the sum on $u$ is over
all the positive integers $u$ such that
$\sqrt{4n+d k^2 u^2}\in\Bbb Z$ and the product on $p^{2l}$ runs over all
distinct primes $p$ with $p^{2l}$ being the greatest even
$p$-power factor of $(d,{N\over k})$, then $L_n(s)$ is analytic for
$\text{Re}\, s>1$ and is extended by analytic continuation to the
half-plane $\text{Re}\, s> 0$ except for a possible pole at $s=1/2$ and for
possible simple poles at $s=1, \frac 12\pm i\kappa_j$,
$j=1,2,\cdots$.  Moreover, we have
$$\text{tr}_{\lambda_j} T_n=2n^{i\kappa_j}
\text{Res}_{s=1/2+i\kappa_j} L_n(s)$$
for $j=1,2,\cdots$.  \endproclaim

  The paper is organized as follows.  In section 2, we recall
the Selberg trace formula; see Hejhal [4] [5] and Selberg [12].
Then in section 3 we
compute contributions of identity, elliptic, hyperbolic
and parabolic elements to the trace formula.   A technical
part of this section is to compute the
Petersson inner product of an Eisenstein series for
$\Gamma_0(N)$ with the image of the Eisenstein series
under the action of Hecke operators.
The result is given in Theorem 3.11.
By using the Selberg trace formula, we obtain
analyticity information about a series formed by contributions
to the trace formula of hyperbolic elements whose fixed points
are not cusps.  A precise statement is given in Theorem 3.13.
In section 4, we compute explicitly the total contribution to
the trace formula of hyperbolic elements whose fixed points
are not cusps.  The result is stated in Theorem 4.6.
Finally, Theorem 1 follows from Theorem 3.13 and Theorem 4.6.

  Most part of this paper was written at the American Institute of
Mathematics, California.   The author wishes to thank AIM and
Brian Conrey for the support.  He also wishes to thank
 Dennis Hejhal and Henryk
Iwaniec for helpful conversations and for their encouragement
for the proof of Theorem 3.11.

\heading
2.   The Selberg trace formula
\endheading

 Assume that $s$ is a complex number with $\text{Re}\, s>1$.  Let
$$k(t)=(1+\frac t4)^{-s}$$
and
$$k(z, z^\prime)=k\left(\frac{|z-z^\prime|^2}{yy^\prime}\right),$$
for $z=x+iy$ and $z^\prime=x^\prime+iy^\prime$ in the upper half-plane.
The kernel function $k(z, z^\prime)$
is of (a)-(b) type in the sense of Selberg [12], p.60.  Let
  $$ g(u)=\int_w^\infty k(t)\frac{dt}{\sqrt {t-w}}$$
with $w=e^u+e^{-u}-2$.   If
  $$h(r)=\int_{-\infty}^{\infty} g(u) e^{iru}du, $$
then
$$g(u)=c(1+\frac w4)^{\frac 12-s} \tag 2.1$$
 where $c=2\sqrt \pi\Gamma(s-\frac 12)\Gamma^{-1}(s)$.   Since
 $$h(r) =c{4^s(s-{1\over 2})\over (s-{1\over 2})^2+r^2}+A(r, s),$$
where $A(r,s)$ is finite for $|\Im r|\leqslant 1/2$ and for
$\text{Re}\, s>0$, we find that
$$\lim_{s\to 1/2+i\kappa} (s-{1\over 2}-i\kappa)h(r)=
\cases 4^{1/2+i\kappa}\sqrt \pi\frac {\Gamma(i\kappa)}{\Gamma(1/2+i\kappa)},
&\text{for $r=\pm\kappa$}; \\0, &\text{for $r\neq \pm\kappa$.}
\endcases\tag 2.2$$

   Assume that $n$ is a positive integer with $(n, N)=1$.  Let
$$\Gamma^*=\cup_{\underset{0\leq b<d}\to{ad=n}}
\frac 1{\sqrt n}\left(\matrix d&{-b}\\0&a\endmatrix\right)\Gamma_0(N).$$
Then $\Gamma^*$ satisfies all the requirements given in [12], p.69.
 Let $\frak a_1, \frak a_2, \cdots, \frak a_{\nu(N)}$
be a complete set of inequivalent cusps of $\Gamma_0(N)$.
We choose $\sigma_{\frak a_i} \in PSL_2(\Bbb R)$  such that
$\sigma_{\frak a_i}\infty=\frak a_i$ and
$\sigma_{\frak a_i}^{-1}\Gamma_{\frak a_i}
\sigma_{\frak a_i}=\Gamma_\infty$ for $i=1,2, \cdots, \nu(N)$.
Let $E_i(z, s)$  be the Eisenstein series for the cusp $\frak a_i$, let
$$K(z, z^\prime)=\sum_{T\in \Gamma^*}k(z, T z^\prime ), $$
and let
$$H(z, z^\prime)= \sum_{i=1}^{\nu(N)}\sum_{ad=n, \,0\leq b< d}
\frac 1{4\pi}\int_{-\infty}^\infty h(r)E_i(\frac{az+b}d, \frac 12+ir)
E_i(z^\prime, \frac 12 -ir) dr.$$
It follows from (2.14) of [12],
the argument of [7], pp.96-98, Theorem 5.3.3 of [7], and the spectral
decomposition formula (5.3.12) of [7] that
$$d(n)h(-\frac i2)+\sqrt n \sum_{j=1}^\infty h(\kappa_j)
\text{tr}_{\lambda_j}T_n=\int_D \{K(z, z)-H(z, z)\}dz \tag 2.3$$
for $\text{Re}\, s>1$, where $d(n)$ is the sum of positive divisors of $n$.

\heading
3.   Evaluation of components of the trace
\endheading

 For every element $T$ of $\Gamma^*$, we denote by $\Gamma_T$ the set
of all the elements of $\Gamma_0(N)$ which commute
with $T$.  Let $D_T=\Gamma_T\backslash \Cal H$.
Elements of $\Gamma^*$ can be classified as
the identity, the hyperbolic, the elliptic, and the parabolic
elements.   If $T$ is not a parabolic element,  we put
$$c(T)=\int_{D_T} k(z, Tz)dz.$$

  Proofs for the following lemmas from Lemma 3.1 to Lemma 3.8
are minor modifications of corresponding lemmas
 in Conrey-Li [2].  For the convenience of readers,
we sketch the proofs here.

\subheading{ 3.1.   The identity component}

  If $\Gamma^*$ contains the identity element $I$, then
  $$c(I)=\int_{\Gamma_0(N)\backslash\Cal H}dz.$$

\subheading{ 3.2.   Elliptic components}

     There are only a finite number of elliptic conjugacy classes.

\proclaim{Lemma 3.1} Let $R$ be an elliptic element of $\Gamma^*$.  Then
$$c(R)=\frac \pi{2m\sin \theta}\int_0^\infty
\frac{k(t)}{\sqrt{t+4\sin^2\theta}}dt, $$
where $m$ is the order of a primitive element of $\Gamma_R$ and
where $\theta$ is defined by the formula trace$(R)=2\cos\theta$.
\endproclaim

  \demo{Proof}   Since $R$ is an elliptic element of $\Gamma^*$,
an element $\sigma\in PSL(2, \Bbb R)$ exists such that
$$\sigma R\sigma^{-1}=\left(\matrix {\cos\theta}&{-\sin\theta}\\
{\sin\theta}&{\cos\theta}\endmatrix\right)=\widetilde R$$
for some number $0<\theta<\pi$.
If $(\sigma\Gamma_0(N)\sigma^{-1})_{\widetilde R}$ is the set
of all the elements of $\sigma\Gamma_0(N)\sigma^{-1}$ commuting
with $\widetilde R$, then we have
$$c(R)=\int_{D_{\widetilde R}}k(z, \widetilde R z)dz$$
where $D_{\widetilde R}
=(\sigma\Gamma_0(N)\sigma^{-1})_{\widetilde R}\backslash \Cal H$.

  If $\gamma$ is an element of $\Gamma_0(N)$ having the same
fixed points as $R$, then $\gamma$ commutes with $R$.
By Proposition 1.16 of [13],
a primitive elliptic element $\gamma_0$ of $\Gamma_0(N)$ exists such that
$(\eta\Gamma\eta^{-1})_{\widetilde R}$ is generated by $\eta\gamma_0\eta^{-1}$.
Since $\eta\gamma_0\eta^{-1}$ commutes with $\widetilde R$, it is of the form
$$\left(\matrix {\cos\theta_0}&{-\sin\theta_0}\\
{\sin\theta_0}&{\cos\theta_0}\endmatrix\right)$$
for some real number $\theta_0$.   By Proposition 1.16 of [13],
$\theta_0=\pi/m$ for some positive integer $m$.
 It follows from the argument of [7], p.99 that
$$c(R)=\frac 1m\int_0^\infty\int_{-\infty}^\infty
k\left(\frac{|z^2+1|^2}{y^2}\sin ^2\theta\right)dz.$$
By the argument of [7], p.100 we have
$$c(R)=\frac \pi{2m\sin \theta}\int_0^\infty
\frac{k(t)}{\sqrt{t+4\sin^2\theta}}dt. \qed$$
\enddemo

\subheading{ 3.3.   Hyperbolic components}

   Let $P$ be a hyperbolic element of $\Gamma^*$.  Then an element $\rho$
exists in $SL_2(\Bbb R)$ such that
$$\rho P\rho^{-1}=\left(\matrix {\lambda_P} &0\\
0&{\lambda_P^{-1}}\endmatrix\right)=\widetilde P$$
with $\lambda_P>1$.  The number $\lambda_P^2$ is called the norm of
$P$, and is denoted by $NP$.  It follows that
$$c(P)=\int_{D_{\widetilde P}}k(z, NP z)dz$$
where $D_{\widetilde P}
=(\rho\Gamma_0(N)\rho^{-1})_{\widetilde P}\backslash \Cal H$.
Let $P_0$ be a primitive hyperbolic element of $SL_2(\Bbb Z)$,
which generates the group of all elements of $SL_2(\Bbb Z)$
commutating with $P$.  Then there exists a ``primitive" hyperbolic
element $P_1\in\Gamma_0(N)$, which generates $\Gamma_P$, such that
$P_1$ is the smallest positive integer power of $P_0$
among all the generators of $\Gamma_P$ in $\Gamma_0(N)$.

\proclaim{Lemma 3.2}  Let $P$ be a hyperbolic element of $\Gamma^*$
such that $\Gamma_P\neq \{1_2\}$.  If $P_1$ is a ``primitive'' hyperbolic
element of $\Gamma_0(N)$ which generates the group $\Gamma_P$, then
$$c(P)=\frac{\ln NP_1}{(NP)^{1/2}-(NP)^{-1/2}}g(\ln NP). $$
\endproclaim

\demo{Proof}   Since $\rho P_1\rho^{-1}$
commutes with $\widetilde P$, it is of the form
$$\left(\matrix {\lambda_{P_1}}&0\\0&{\lambda_{P_1}^{-1}}\endmatrix\right)$$
for some number $\lambda_{P_1}>1$, and hence
$$c(P)=\int_1^{NP_1}\frac{dy}{y^2}\int_{-\infty}^\infty
k\left(\frac{(NP-1)^2}{NP}\frac{|z|^2}{y^2}\right)dx.$$
The stated identity then follows.  \qed\enddemo

\proclaim{Lemma 3.3}  Let $P={1\over \sqrt n}\left(\smallmatrix
A&B\\C&D\endsmallmatrix\right)$ be a hyperbolic element of $\Gamma^*$
such that $\Gamma_P=\{1_2\}$.  Then fixed points of $P$ are cusps
of $\Gamma_0(N)$.  Moreover, $\Gamma_P=\{1_2\}$
if and only if $A+D={1\over 2}(m+{4n\over m})$ for some divisor $m$ of
$4n$ with $m\neq 2\sqrt n$ for $C\neq 0$, and if and only if $A\neq D$
for $C=0$.
\endproclaim

 \demo{Proof}  If the fixed points of $P$ are not
rational numbers or infinity, then they are zeros of an irreducible polynomial
$ax^2+bx+c$ with $a\equiv 0$ mod$(N)$.  If $d=b^2-4ac$ and $(u, v)$ is
a solution of Pell's equation $v^2-du^2=4$, then
$$\left(\matrix {v-bu\over 2}&{-cu}\\{au}&{v+bu\over 2}\endmatrix\right)$$
belongs to $\Gamma_0(N)$ and commutes with $P$.
This contradicts to $\Gamma_P=\{1_2\}$.
Hence, fixed points of $P$ are cusps of $\Gamma_0(N)$.

  Conversely, if $T$ is an element of $\Gamma^*$ having two distinct
fixed points with at least one of them being a rational number,
then we must have $\Gamma_T=\{1_2\}$.   Hence,
$\Gamma_P=\{1_2\}$ for an element $P={1\over \sqrt n}
\left(\smallmatrix A&B\\C&D\endsmallmatrix\right)\in\Gamma^*$
if and only if the fixed points of $P$ are rationals or the infinity,
that is, if and only if $A+D={1\over 2}(m+{4n\over m})$ for some
divisor $m$ of $4n$ with $m\neq 2\sqrt n$ if $C\neq 0$ or $A\neq D$
if $C=0$.  \qed\enddemo

 Every cusp of $\Gamma_0(N)$ is equivalent to one of
the following inequivalent cusps
$$ {u\over w}\,\,\text{with}\,\, u, w>0,
\,\, (u, w)=1,\,\, w|N.\,\,\tag 3.1$$
Two such cusps $u/w$ and $u_1/w_1$ are $\Gamma_0(N)$-equivalent
if and only if $w=w_1$ and $u\equiv u_1$ modulo $(w, N/w)$.
Let $\frak a=u/w$ be given as in (3.1).  By (2.2) and (2.3)
of Deshouillers and Iwaniec [3], we have
$$\Gamma_{\frak a}=\left\{\left(\matrix {1+cu/w}&{-cu^2/w^2}\\c&{1-cu/w}
\endmatrix\right):\,\, c\equiv 0
\,\,(\text{mod}\,[w^2, N])\right\} \tag 3.2$$
and
$$\sigma_{\frak a}\infty=\frak a\,\,\,\,\text{and}\,\,\,\,
\sigma_{\frak a}^{-1}\Gamma_{\frak a}\sigma_{\frak a}=\Gamma_\infty$$
where
$$\sigma_{\frak a}=\left(\matrix {\frak a\sqrt{[w^2, N]}}&0\\
{\sqrt{[w^2, N]}}&{1/\frak a \sqrt{[w^2, N]}}\endmatrix\right).
\tag 3.3$$
Let $P$ be a hyperbolic element of $\Gamma^*$
such that $\Gamma_P=\{1_2\}$.  Assume that
$\frak a$ is a fixed point of $P$.  Then $\infty$ is a fixed point
of $\sigma_{\frak a}^{-1}P\sigma_{\frak a}$, and hence there exist
positive integers $a,d$ with $ad=n, a\neq d$ such that
$$P={1\over\sqrt n}\sigma_{\frak a}\left(\matrix a&b\\0&d\endmatrix\right)
\sigma_{\frak a}^{-1}={1\over \sqrt n}\left(\matrix{a-b[w^2, N]
{u\over w}}&{b[w^2, N]{u^2\over w^2}}\\
{(a-d){w\over u}-b[w^2, N]}&{d+b[w^2, N]{u\over w}}
\endmatrix\right).\tag 3.4$$
Since $c(P)$ depends only on the conjugacy class
$\{P\}$ represented by $P$, we can replace $P$
 by $\gamma^{-1}P\gamma$ without changing the value of $c(P)$.
Replacing $P$ by $\gamma^{-1}P\gamma$ for some element
$\gamma=\sigma_{\frak a}\left(\smallmatrix 1& \ell \\ 0&1\endsmallmatrix
\right)\sigma_{\frak a}^{-1}\in \Gamma_{\frak a}$,
we can assume without loss of generality that
$0\leqslant b<|a-d|$ in (3.4).

  \proclaim{Lemma 3.4}  Let $\frak a=u/w$ be given as in (3.1).  Then
$$P={1\over\sqrt n}\sigma_{\frak a}\left(\matrix a&b\\0&d\endmatrix\right)
\sigma_{\frak a}^{-1} $$
with $a, d\in\Bbb Z, ad=n, a\neq d$ is a hyperbolic element of $\Gamma^*$ with
$\Gamma_P=\{1_2\}$ and $P(\frak a)=\frak a$ if and only if
$(w, N/w)|(a-d)$ with $b$ being chosen so that $P\in\Gamma^*$.\endproclaim

  \demo{Proof}  Let $P$ be a hyperbolic element of $\Gamma^*$ with
$\Gamma_P=\{1_2\}$ and $P(\frak a)=\frak a$.  By (3.4), there exists
an integer $k$ such that ${u\over w}N|a-d-kw$, and hence $(w, N/w)|a-d$.

   Conversely, if $a, d$ are positive integers satisfying $ad=n, a\neq d$
and $(w, N/w)|a-d$, then  $a-d=\ell (w, N/w)$ for some integer $\ell$.
Integers $\lambda$ and $\tau$ exist such that
$\lambda\frak a N+\tau w=(w, N/w)$, and hence
$\ell\lambda \frak a N=a-d-\ell\tau w$.  Thus, an integer $k$ exists
such that ${u\over w}N|a-d-kw$.  Let
$A=a-kw, B=ku, C={w\over u}(a-d-kw)$ and $D=d+kw$.  Then
$P={1\over\sqrt n}\left(\smallmatrix A&B\\C&D\endsmallmatrix\right)$
is an element of $\Gamma^*$ with $\Gamma_P=\{1_2\}$ and $P(\frak a)=\frak a$,
which is of the form
$$P={1\over\sqrt n}\sigma_{\frak a}\left(\matrix a&b\\0&d\endmatrix\right)
\sigma_{\frak a}^{-1} $$
with $a, d\in\Bbb Z, ad=n, a\neq d$ and with $b$ being chosen so that
$P\in\Gamma^*$.  \qed\enddemo

  \proclaim{Lemma 3.5}   Let $\frak a=u/w$ be given as in (3.1).
Let $a, d\in\Bbb Z^+, ad=n, a\neq d$.  If $(w, N/w)|a-d$, then there
are exactly $|a-d|$ number of $\Gamma_0(N)$-inequivalent hyperbolic
elements $P\in\Gamma^*$ with $P(\frak a)=\frak a$ and $\Gamma_P
=\{1_2\}$, which are of the form
$$P={1\over\sqrt n}\sigma_{\frak a}\left(\matrix a&b\\0&d\endmatrix\right)
\sigma_{\frak a}^{-1}$$
with $0\leqslant b<|a-d|$.  \endproclaim

  \demo{Proof}    For $a, d\in\Bbb Z, ad=n, a\neq d$, let
$$P={1\over\sqrt n}\sigma_{\frak a}\left(\matrix a&b\\0&d\endmatrix\right)
\sigma_{\frak a}^{-1}\,\,\,\,\text{and}\,\,\,\,
P^\prime={1\over\sqrt n}\sigma_{\frak a}\left(\matrix a&b^\prime\\0&d
\endmatrix\right) \sigma_{\frak a}^{-1}$$
be two elements of $\Gamma^*$ with $0\leqslant b, b^\prime<|a-d|$
whose fixed points are cusps.  By using (3.4), we can find that
$b^\prime=b+\ell$ for some integer $\ell$.

 Conversely, for a given integer $\ell$ with
$0\leqslant b+\ell< |a-d|$, if $b^\prime=b+\ell$ then
$$P^\prime={1\over\sqrt n}\sigma_{\frak a}\left(\matrix
a&b^\prime\\0&d\endmatrix\right) \sigma_{\frak a}^{-1}$$
is a hyperbolic element of $\Gamma^*$ with
$\Gamma_{P^\prime}=\{1_2\}$ and $P^\prime(\frak a)=\frak a$.
Moreover, if $b^\prime\neq b$ modulo $|a-d|$, then
$P^\prime$ and $P$ are not $\Gamma_0(N)$-equivalent.
Therefore, for fixed $a, d\in\Bbb Z^+, ad=n, a\neq d$, there
are exactly $|a-d|$ $\Gamma_0(N)$-inequivalent hyperbolic
elements $P$ in $\Gamma^*$ with $P(\frak a)=\frak a$ and $\Gamma_P
=\{1_2\}$, which are of the form
$$P={1\over\sqrt n}\sigma_{\frak a}\left(\matrix a&b\\0&d\endmatrix\right)
\sigma_{\frak a}^{-1}$$
with $0\leqslant b<|a-d|$.  \qed\enddemo

  \proclaim{Lemma 3.6}  Let $\frak a=u/w$, and let
$$P={1\over\sqrt n}\sigma_{\frak a}\left(\matrix a&b\\0&d\endmatrix\right)
\sigma_{\frak a}^{-1}\,\,\,\,\text{and}\,\,\,\,
P^\prime={1\over\sqrt n}\sigma_{\frak a}\left(\matrix d&b^\prime\\0&a
\endmatrix\right) \sigma_{\frak a}^{-1}$$
be two hyperbolic elements of $\Gamma^*$ with  $P(\frak a)=\frak a,
P^\prime(\frak a)=\frak a$, $\Gamma_P=\{1_2\}$ and $\Gamma_{P^\prime}
=\{1_2\}$.   Then $P$ is $\Gamma_0(N)$-conjugate to $P^\prime$ for some
number $b^\prime$ if and only if the two fixed points of $P$
are $\Gamma_0(N)$-conjugate.   \endproclaim

 \demo{Proof}   Assume that the two fixed points of $P$
are $\Gamma_0(N)$-conjugate.  Then there exists an
element $T=\left(\smallmatrix \alpha&\beta\\
\gamma&\delta\endsmallmatrix \right)\in\Gamma_0(N)$
such that $\frak a, T(\frak a)$ are the two
fixed points of $P$.   Let $P^\prime=T^{-1}PT$.  Then
$P^\prime(\frak a)=\frak a$ and $\Gamma_{P^\prime}=\{1_2\}$.
  Write
$$P^\prime={1\over\sqrt n}\sigma_{\frak a}\left(\matrix
{a^\prime}&{b^\prime}\\0&{d^\prime}\endmatrix\right)
\sigma_{\frak a}^{-1}$$
with $a^\prime, d^\prime\in\Bbb Z^+, a^\prime d^\prime=n$, and
$0\leqslant b^\prime<|a^\prime-d^\prime|$.
Since $P^\prime=T^{-1}PT$, we have
$$(\sigma_{\frak a}^{-1}T\sigma_{\frak a})\left(\matrix
a^\prime&b^\prime\\0&d^\prime\endmatrix\right)
=\left(\matrix a&b\\0&d\endmatrix\right)
(\sigma_{\frak a}^{-1}T\sigma_{\frak a}). $$
From this identity, we can deduce that $a^\prime=d$, and hence
$d^\prime=a$.

  Conversely, suppose that there exists an element
 $T=\left(\smallmatrix A&B\\C&D\endsmallmatrix\right)\in
\Gamma_0(N)$ such that $T^{-1}PT=P^\prime$, that is,
$$\sigma_{\frak a}^{-1}T^{-1}\sigma_{\frak a}
\left(\matrix a&b\\0&d\endmatrix\right)
=\left(\matrix d&{b^\prime}\\0&a\endmatrix\right)
\sigma_{\frak a}^{-1}T^{-1}\sigma_{\frak a}.$$
It follows that
$$\left((a-d){w\over u}-b[w^2, N]\right)u(A+B
{w\over u})=b[w^2, N]{u\over w}(-C u-D w).$$
This identity can be written as
$$T(\frak a)={\frak a\over 1-{a-d\over b[w^2, N]\frak a}},$$
which is the second fixed point of $P$.  That is, the two fixed
points of $P$ are $\Gamma_0(N)$-conjugate.  \qed\enddemo

\proclaim{Lemma 3.7}  Let $a,d\in \Bbb Z^+, a\neq d, ad=n$,
and let $\frak a=u/w$ be given as in (3.1) with
$(w, N/w)|(a-d)$.  Assume that $\frak a$ and
$\frak a^\prime=u_1/w_1$ are the two distinct fixed point of
$P={1\over\sqrt n}\sigma_{\frak a}\left(\smallmatrix a&b
\\0&d\endsmallmatrix\right)\sigma_{\frak a}^{-1}$.
Then $(w^\prime, N/w^\prime)|(a-d)$ for $w^\prime=(w_1, N)$, and
$$P={1\over\sqrt n}\sigma_{\frak a^\prime}\left(\matrix d&b^\prime
\\0&a\endmatrix\right)\sigma_{\frak a^\prime}^{-1}$$
for some number $b^\prime$. \endproclaim

  \demo{Proof}  There exist positive
integers $a^\prime, d^\prime$ with $a^\prime d^\prime=n$ such that
$$P={1\over\sqrt n}\sigma_{\frak a^\prime}\left(\matrix a^\prime
&b^\prime \\0&d^\prime\endmatrix\right)\sigma_{\frak a^\prime}^{-1}$$
for some number $b^\prime$, and hence we have
$$\sigma_{\frak a}^{-1}\sigma_{\frak a^\prime}\left(\matrix a^\prime
&b^\prime \\0&d^\prime\endmatrix\right)
=\left(\matrix a&b\\0&d\endmatrix\right)
\sigma_{\frak a}^{-1}\sigma_{\frak a^\prime}. \tag 3.5$$
Write $\sigma_{\frak a}^{-1}\sigma_{\frak a^\prime}=\left(
\smallmatrix A&B\\C&D\endsmallmatrix\right)$.  Then by
(3.5) we get that $(d-a^\prime)C=0, (a^\prime -a)A=bC$
and $(d-d^\prime)D=b^\prime C$.  Since $\frak a\neq\frak a^\prime$,
by using (3.3) we find that $C\neq 0$, and hence we have
$a^\prime=d$.  It follows that $d^\prime=a$.
By the argument made in the first paragraph of the
proof of Lemma 3.4, we have $(w^\prime, N/w^\prime)|(a-d)$.
 \qed\enddemo

   For a large positive number $Y$, define
$$D_Y=\{z\in D: \Im \sigma_{\frak a_i}^{-1}z<Y,\,\, i=1,2,\cdots,h\}.$$
Let
$$(D_P)_Y=\bigcup_{\gamma\in\Gamma_0(N)} \gamma D_Y.$$
Write
$$c(P)_Y=\int_{(D_P)_Y}k(z, Pz)dz. $$

 \proclaim{Lemma 3.8}  Let $\frak a=u/w$ be given as in (3.1), and let
$P={1\over\sqrt n}\sigma_{\frak a}\left(\smallmatrix a&b
\\0&d\endsmallmatrix\right)\sigma_{\frak a}^{-1}$
with $a, d\in\Bbb Z^+, ad=n, a\neq d$ be a hyperbolic element
of $\Gamma^*$ with $\Gamma_P=\{1_2\}$ and $P(\frak a)=\frak a$.
Then
$$c(P)_Y=\frac{\sqrt n\ln \{(a-d)^2[w^2, N]Y/2\rho\}}{|a-d|}
g(\ln \frac ad) +\int_1^\infty  k\left(\frac{(a-d)^2}nt\right)
\frac{\ln t}{\sqrt {t-1}}dt+o(1) $$
where  $o(1)\to 0$ as $Y\to\infty$ and
$\rho=C^2/2[C^2/\ell^2, N]Y$ with $C=(a-d){w\over u}-b[w^2, N]$
and $\ell=\text{gcd}\left(C, b[w^2,N]{u\over w}\right)$.  \endproclaim

  \demo{Proof}   Let
$$\gamma=\left(\matrix p&r\\q&s\endmatrix \right)$$
be an element of $SL_2(\Bbb R)$.  The linear fractional
transformation, which takes
every complex $z$ in the upper half-plane into $\gamma (z)$, maps the
horizontal line $\Im z=Y$ into a circle of radius $\frac 1{2q^2 Y}$
with center at $\frac pq+\frac i{2q^2 Y}$.
Let $$\mu=\left(\matrix 1&{\frac b{a-d}}\\0&1\endmatrix \right).$$
Then
$$P={1\over\sqrt n}\sigma_{\frak a}\mu^{-1}\left(\matrix a&0\\0&d
\endmatrix\right)\mu\sigma_{\frak a}^{-1}.$$
Note that
$$\mu\sigma_{\frak a}^{-1}=\left(\matrix {{w/u\over\sqrt{[w^2, N]}}
-{b\sqrt{[w^2, N]}\over a-d}}&{bu\sqrt{[w^2, N]}\over (a-d)w}
\\{-\sqrt{[w^2, N]}}&{{u\over w}\sqrt{[w^2, N]}}\endmatrix\right).$$
Since
$$P={1\over\sqrt n}\sigma_{\frak a}\mu^{-1}\left(\matrix {1/y}&0\\
0&y\endmatrix\right)\left(\matrix a&0\\0&d \endmatrix\right)
\left(\matrix y&0\\0&{1/y}\endmatrix\right)\mu\sigma_{\frak a}^{-1}$$
for any $y\neq 0$, by choosing $y=(a-d)\sqrt{[w^2, N]}$ we obtain that
 $$c(P)_Y=\int_{\mu_{\frak a}\{(D_P)_Y\}} k(z, {a\over d}z)dz$$
where
$$\mu_{\frak a}=\left(\matrix {(a-d){w\over u}-b[w^2, N]}
&{b[w^2, N]{u\over w}}\\{-1\over a-d}&{u\over (a-d)w}
\endmatrix\right).$$

The linear transformation $z\to\sigma_{\frak a} z$ maps the
half-plane $\Im z>Y$ into a disk $D_{\frak a}$ of radius
${1\over 2[w^2, N]Y}$ with center at $\frak a+{i\over 2[w^2, N]Y}$.
Then the transformation $z\to \mu_{\frak a}(z)$
maps the disk $D_{\frak a}$ into the half-plane
$\Im z> (a-d)^2[w^2, N]Y$.
Let
$$\frak p={-b[w^2, N]u/w\over (a-d){w\over u}-b[w^2, N]}.$$
If
$$\sigma_{\frak p}=\left(\matrix {\frak p\sqrt
{[C^2/\ell^2, N]}}&0\\ {\sqrt{[C^2/\ell^2, N]}}&
{1/\frak p\sqrt{[C^2/\ell^2, N]}}\endmatrix\right),$$
 then $\sigma_{\frak p}\infty=\frak p$
and $\sigma_{\frak p}^{-1}\Gamma_{\frak p}\sigma_{\frak p}
=\Gamma_{\infty}$.  By the definition of $D_Y$,
 the image of the half-plane $\Im z>Y$ under the
linear transformation $z\to\sigma_{\frak p}z$ is not
contained in $(D_P)_Y$.   Since the linear transformation
$z\to (\mu_{\frak a}\sigma_{\frak p})(z)$ maps the
half-plane $\Im z>Y$ into a disk $D_\rho$ of radius
$\rho$ centered at $i\rho$, where $\rho=C^2/2[C^2/\ell^2, N]Y$,
the disk $D_\rho$ is not contained in
$\mu_{\frak a}\{(D_P)_Y\}$.  It follows that
$$\aligned c(P)_Y &=\int_0^\pi d\theta
\int_{2\rho \sin \theta}^{(a-d)^2[w^2, N]Y/\sin \theta}k
\left({(a-d)^2\over n\sin^2\theta}
\right)\frac{dr}{r\sin^2\theta}+o(1)\\
&=\int_1^\infty k\left(\frac{(a-d)^2}nt\right)
\frac{\ln \{(a-d)^2[w^2, N]Yt/2\rho\}}{\sqrt {t-1}}dt
+o(1) \endaligned $$
where $o(1)$ has a limit zero when $Y\to\infty$, and hence we have
$$c(P)_Y=\frac{\sqrt n\ln \{(a-d)^2[w^2, N]Y/2\rho\}}{|a-d|}
g(\ln \frac ad) +\int_1^\infty  k\left(\frac{(a-d)^2}nt\right)
\frac{\ln t}{\sqrt {t-1}}dt+o(1). $$
This completes the proof of the lemma. \qed\enddemo

\proclaim{Theorem 3.9}  Let $N$ be any positive integer, and let
$$H=\{(\frak a; a,b,d): \text{$\frak a, a,b$, and $d$ are given
as in Lemma 3.4}\}.$$
  Then we have
$$\aligned  \sum_{\{P\}, \Gamma_P=\{1_2\}} &c(P)_Y
=\sqrt n\ln Y\sum_{\underset{(w, N/w)|(a-d)}\to{{u\over w},
ad=n, d>0, a\neq d}}g(\ln{a\over d})\\
&+{\sqrt n\over 2}\sum_{(\frak a; a,b,d)\in H}
\ln \{(a-d)^2[w^2,N][{C^2\over\ell^2},N]/C^2\}
g(\ln \frac ad)\\
&+{1\over 2}\sum_{\underset{(w, N/w)|(a-d)}\to{{u\over w},
ad=n, d>0, a\neq d}} |a-d|\int_1^\infty  k\left(\frac{(a-d)^2}nt\right)
\frac{\ln t}{\sqrt {t-1}}dt+o(1)\endaligned$$
where $C=(a-d){w\over u}-b[w^2, N]$ and
$\ell=\text{gcd}\left(C, b[w^2,N]{u\over w}\right)$.  \endproclaim

  \demo{Proof}  Let $S_1$ be the set of conjugacy classes of
hyperbolic elements $P$ of $\Gamma^*$ such that $\Gamma_P=\{1_2\}$
and such that the two distinct cusps of $P$ are not
$\Gamma_0(N)$-equivalent, and let $S_2$ be the set of conjugacy
classes of hyperbolic elements $P$ of $\Gamma^*$
such that $\Gamma_P=\{1_2\}$ and such that the two distinct
cusps of $P$ are $\Gamma_0(N)$-equivalent.  For the convenience
of writing, we denote
$$P(\frak a; a, b, d)={1\over\sqrt n}\sigma_{\frak a}
\left(\matrix a&b\\0&d\endmatrix\right)\sigma_{\frak a}^{-1}$$
where $\frak a=u/w$ is given as in (3.1).  Since the two
distinct cusps of every element $P$ in $S_1$ are
equivalent to $\Gamma_0(N)$-inequivalent cusps of the form (3.1),
we have
$$\sum_{\{P\}\in S_1} c(P)_Y
={1\over 2}\sum_{\underset{P(\frak a; a, b, d)\in S_1}\to
{(\frak a; a,b,d)\in H}}c(P(\frak a; a, b, d))_Y$$
by Lemma 3.6.  It follows from Lemma 3.6 and Lemma 3.7 that
$$\sum_{\{P\}\in S_2} c(P)_Y
={1\over 2}\sum_{\underset{P(\frak a; a, b, d)\in S_2}\to
{(\frak a; a,b,d)\in H}}c(P(\frak a; a, b, d))_Y.$$
It follows that
$$\sum_{(\frak a; a,b,d)\in H}c(P(\frak a; a, b, d))_Y
=2\sum_{\{P\}, \Gamma_P=\{1_2\}} c(P)_Y.$$
By Lemma 3.5 and Lemma 3.8, we have
$$\aligned &\sum_{(\frak a; a,b,d)\in H}c(P(\frak a; a, b, d))_Y
=2\sqrt n\ln Y\sum_{\underset{(w, N/w)|(a-d)}\to{{u\over w},
ad=n, d>0, a\neq d}}g(\ln{a\over d})\\
&+\sqrt n\sum_{(\frak a; a,b,d)\in H}
\ln \{(a-d)^2[w^2,N][{C^2\over\ell^2},N]/C^2\} g(\ln \frac ad)\\
&+\sum_{\underset{(w, N/w)|(a-d)}\to{{u\over w},
ad=n, d>0, a\neq d}} |a-d|\int_1^\infty  k\left(\frac{(a-d)^2}nt\right)
\frac{\ln t}{\sqrt {t-1}}dt+o(1). \qed\endaligned$$
 \enddemo

\subheading{ 3.4.   Parabolic components}

  Let $S$ be a parabolic element of $\Gamma^*$.   Then
every element of $\Gamma_0(N)$ having the same fixed point as
$S$ commutes with $S$.   If $\frak a=u/w$ is the fixed point
of $S$, then $\Gamma_S=\Gamma_{\frak a}$, and hence we have
$\sigma_{\frak a}^{-1}\Gamma_S \sigma_{\frak a}=\Gamma_\infty$,
where $\sigma_{\frak a}$ is given as in (3.3).
Since $\sigma_{\frak a}^{-1} S\sigma_{\frak a}$ commutes
with every element of $\Gamma_\infty$, we have
$$\sigma_{\frak a}^{-1} S\sigma_{\frak a}={1\over \sqrt n}
\left(\matrix a&b\\0&a\endmatrix\right)$$
for some numbers $a, b$ with $a^2=n$, $a\in\Bbb Z$.
This implies that  $\Gamma^*$ has  parabolic elements only
if $n$ is the square of an integer.
Furthermore,  by (3.4) we see that elements of the form
$$S=\sigma_{\frak a}\left(\matrix 1&{b/{\sqrt n}}\\0&1\endmatrix\right)
\sigma_{\frak a}^{-1}, \, \, \, \, \, \, \,0\neq b\in \Bbb Z$$
constitute a complete set of representatives for the conjugacy
classes of parabolic elements of $\Gamma^*$ having $\frak a$ as
its fixed point.  It follows that
$$\sum_{\{S\}}\int_{D_Y}k(z, Sz)dz =\int_0^Y\int_0^1
\sum_{0\neq b\in \Bbb Z}k(z, z+\frac b{\sqrt n})dz+o(1), $$
where the summation on $\{S\}$ is taken over all parabolic
classes represented by parabolic elements whose fixed point is
$\frak a=u/w$ and where $o(1)$ tends to zero as $Y\to\infty$.
Define $\delta_n$ to be one if $n$ is the square of an
integer and to be zero otherwise.

   The identity
$$\sigma_{\frak a_i}\left(\matrix A&B\\0&D
\endmatrix\right)\sigma_{\frak a_j}^{-1}
=\left(\matrix {A{\frak a_i\over\frak a_j}-B\frak a_i[w^2, N]}
&{B\frak a_i\frak a_j[w^2, N]}\\
{{A\over\frak a_j}-{D\over\frak a_i}-B[w^2,N]}
&{\frak a_j B[w^2,N]+D{\frak a_j\over\frak a_i}}\endmatrix\right)
\tag 3.6$$
holds for $\frak a_i=u_i/w$ and $\frak a_j=u_j/w$.

\proclaim{Lemma 3.10} If
$${1\over\sqrt n}\left(\matrix A&B\\0&D\endmatrix\right)$$
is an element of
$\sigma_{\frak a_i}^{-1}\Gamma^*\sigma_{\frak a_i}$,
 then $A, D$ are integers satisfying
 $(w, {N\over w})|(A-D)$.  Conversely, if $(w, {N\over w})|(A-D)$,
 then there are exactly $D$ elements in $\Gamma_{\frak a_i}\backslash\Gamma^*$
 represented by matrices of the form
${1\over\sqrt n}\sigma_{\frak a_i}\left(
\smallmatrix A&B\\0&D\endsmallmatrix\right)\frak a_i^{-1}$.
 \endproclaim

\demo{Proof}  Assume that an element $T=\left(\smallmatrix \alpha&\beta \\
\gamma&\delta\endsmallmatrix\right)\in\Gamma_0(N)$ exists such that
 $$\sigma_{\frak a_i}^{-1}T\left(\matrix a&b\\0&d\endmatrix\right)
\sigma_{\frak a_i}= \left(\matrix A&B\\0&D\endmatrix\right)$$
for some positive integers $a,d$ with $ad=n$.  Then
$$B={(\alpha b+\beta d)w^2\over u_i^2[w^2,N]}.$$
Since
$${1\over\sqrt n}\sigma_{\frak a_i}\left(\matrix A&B\\0&D
\endmatrix\right)\sigma_{\frak a_i}^{-1}$$
is an element of $\Gamma^*$, $A,D$ must be integers.
By (3.6) we have that

$$B\frak a_i^2[w^2,N]={Bu_i^2[w^2,N]\over w^2}$$
is an integer.  Since $(u_i, w)=1$, it follows
that
$${B[w^2,N]\over w}\equiv 0\,\, (\text{mod} \, w).$$
 By (3.6) we also have that $N/w$ divides
$${A-D\over u_i}-{B[w^2,N]\over w}.$$
It follows that
$${A-D\over u_i}\equiv 0\,\, (\text{mod}\,
(w, {N\over w})).$$
Since $(u_i, w)=1$, we have that
$$A-D\equiv 0\,\, (\text{mod}\,(w, {N\over w})).$$

  Conversely, if $(w, {N\over w})|(A-D)$, then it follows from
the second part of the proof of Lemma 3.4 that
$${1\over\sqrt n}\sigma_{\frak a_i}
\left(\matrix A&B\\0&D\endmatrix\right)
\sigma_{\frak a_i}^{-1}$$
is an element of $\Gamma^*$ for some number of $B$.
By an argument similar to that made in the proof of the lemma 3.5,
there are exactly $D$ elements in $\Gamma_{\frak a_i}\backslash\Gamma^*$
represented by matrices of the form
${1\over\sqrt n}\sigma_{\frak a_i}\left(
\smallmatrix A&B\\0&D\endsmallmatrix\right)\frak a_i^{-1}$.

  This completes the proof of the lemma. \qed\enddemo

\proclaim{The Maass-Selberg relation [7]}
Let $\Gamma$ be a discontinuous group of finite type, and
let $\frak a_1, \cdots, \frak a_h$ be a complete set of inequivalent
cusps of $\Gamma$.  Assume that $f(z, s)$ and $g(z, s^\prime)$
are two automorphic functions with respect to $\Gamma$ such that
$$\Delta f(z, s)=s(1-s)f(z,s)$$
and
$$\Delta g(z, s^\prime)=s^\prime(1-s^\prime) g(z, s^\prime)$$
and such that the constant terms of the Fourier expansion
at $\frak a_j$ of $f$ and $g$ are given
respectively by
$$c_jy^s+c^\prime_j y^{1-s}$$
and
$$d_jy^{s^\prime}+d^\prime_j y^{1-s^\prime}.$$
Then, for the compact parts $f^Y$ and $g^Y$ of $f$ and $g$,
we have
$$\aligned &\langle f^Y(z, s), \bar g^Y(z, s^\prime)\rangle\\
&=\sum_{j=1}^h\left({c_jd_jY^{s+s^\prime-1}-c_j^\prime
d_j^\prime Y^{-s-s^\prime+1}\over s+s^\prime -1}+
{c_jd_j^\prime Y^{s-s^\prime}-c_j^\prime d_j Y^{-s+s^\prime}
\over s-s^\prime}\right)\endaligned $$
whenever $s(1-s)\neq s^\prime(1-s^\prime)$ and
$f^Y$ and $g^Y$ are square integrable on $\Gamma\backslash
\Cal H$.  \endproclaim

\proclaim{Theorem 3.11}  We have
$$\aligned &\sum_{i=1}^{\nu(N)}\int_{D_Y}
\{\sum_{ad=n, 0\leqslant b<d}E_i\left({az+b\over d},s\right)\}
E_i(z, \bar s)dz \\
&=\sum_{j=1}^{\nu(N)}\sum_{\underset{\frak a_j={u_j\over
w_j}, (w_j, N/w_j)|(a-d)}\to{a,d>0, ad=n}}a^sd^{1-s}
 \{{Y^{s+\bar s-1}-\sum_{i=1}^{\nu(N)}\varphi_{ij}(s)
\varphi_{ij}(\bar s) Y^{-s-\bar s+1}\over s+\bar s-1}\\
&+{\varphi_{jj}(\bar s) Y^{s-\bar s}
-\varphi_{jj}(s) Y^{-s+\bar s}\over s-\bar s}\}
-\sum_{j=1}^{\nu(N)}\sum_{a,d,b}^*
{\left({A^\prime \over D^\prime}\right)^{1-s}\varphi_{jj^\prime}(s)
 Y^{-s+\bar s}\over s-\bar s}+o_Y(1)\endaligned $$
as $\text{Re}\, s\to 1/2$,
where the sum on $a,d,b$ has the same meaning as in (3.10)
with $A^\prime, D^\prime$ being defined there
and $o_Y(1)$ tends to zero as $Y\to\infty$.
\endproclaim

\demo{Proof}  Let
 $$f_i(z)=\sum_{ad=n, 0\leqslant b<d}E_i\left({az+b\over d},s\right).$$
Then $f_i(z)$ is an automorphic function with respect to $\Gamma_0(N)$.
Let the constant term of the Fourier expansion at $\frak a_j$ of
$f_i(z)$ be $c_{ij}y^s+c_{ij}^\prime y^{1-s}$.
By the Maass-Selberg relation we have
$$\aligned &\sum_{i=1}^{\nu(N)}
\int_{D_Y} f_i(z)E_i(z, \bar s)dz-o_Y(1)  \\
&=\sum_{i,j=1}^{\nu(N)}\left({c_{ij}\delta_{ij}
Y^{s+\bar s-1}-c_{ij}^\prime
\varphi_{ij}(\bar s) Y^{-s-\bar s+1}\over s+\bar s-1}+
{c_{ij}\varphi_{ij}(\bar s) Y^{s-\bar s}-c_{ij}^\prime \delta_{ij}
 Y^{-s+\bar s}\over s-\bar s}\right),\endaligned
 \tag 3.7$$
where $o_Y(1)$ tends to zeros as $Y\to\infty$.

Since
$$E_i\left({a\sigma_{\frak a_j}z+b\over d},s\right)
=\sum_{\gamma\in\Gamma_{\frak a_i}\backslash \Gamma_0(N)}
y\left(\sigma_{\frak a_i}^{-1}
\gamma{a\sigma_{\frak a_j}z+b\over d}\right)^s,$$
we have
$$f_i(\sigma_{\frak a_j}z)=\sum_{ad=n, 0\leqslant b<d}
\sum_{\gamma\in\Gamma_{\frak a_i}\backslash \Gamma_0(N)}
y\left(\sigma_{\frak a_i}^{-1}
\gamma\left(\smallmatrix a&b\\0&d\endsmallmatrix\right)
\sigma_{\frak a_j}z\right)^s.$$

  Let $a,d>0$ and $b$ be integers with $ad=n$ and $0\leq b<d$.  If
the matrix equation
$$\sigma_{\frak a_j}^{-1}{1\over\sqrt n}T\left(\matrix a&b\\
0&d\endmatrix\right)\sigma_{\frak a_j}=
{1\over\sqrt n}\left(\matrix A&B\\0&D\endmatrix\right)
\tag 3.8$$
holds for an element $T$ in $\Gamma_0(N)$ and
for some numbers $A, D$ and $B$, then $A,D$ are integers with
$AD=n$.  Moreover, by Lemma 3.10 we must have $(w_i, N/w_i)|(A-D)$ for
$\frak a_i=u_i/w_i$.   Note that, for any fixed pair of positive
integers $A,D$ with $AD=n$, if
the matrix equation (3.8) has a solution $T, a,b,d$ for a
number $B$ with $0\leq B<D$, then the matrix equation (3.8) has solutions
for exactly $D$ number of $B$'s with $0\leq B<D$ by the proof of Lemma 3.5.
By Theorem 2.7 of Iwaniec [6], we have the following decomposition
$$\aligned &\Gamma_\infty\backslash\sigma_{\frak a_i}^{-1}
{1\over\sqrt n}\Gamma_0(N)\left(\matrix
a&b\\0&d\endmatrix\right)\sigma_{\frak a_j}\\
&=\left[\delta_{ij} \cup \bigcup_{\gamma>0}
\bigcup_{\delta\,(\text{mod}\,\gamma)}
\left(\matrix *&*\\ \gamma &\delta
\endmatrix\right)\Gamma_\infty\right]{1\over\sqrt n}
\left(\matrix A&B\\ 0&D\endmatrix\right)\endaligned $$
where the unions on $\gamma, \delta$ run over all numbers $\gamma>0$,
$0\leq \delta<\gamma$ such that
$\left(\smallmatrix *&*\\ \gamma&\delta
\endsmallmatrix\right)$ belong to $\sigma_{\frak a_i}^{-1}
\Gamma_0(N)\sigma_{\frak a_j}$.   Then
it follows from Lemma 3.10 and the proof of
Theorem 3.4 of Iwaniec [6]
that the constant coefficient in the Fourier
expansion at infinity of the function
$$\widetilde f_i(\sigma_{\frak a_j}z)
=\sum_*\sum_{\gamma\in\Gamma_{\frak a_i}
\backslash \Gamma_0(N)}
y\left(\sigma_{\frak a_i}^{-1}
\gamma\left(\smallmatrix a&b\\0&d\endsmallmatrix\right)
\sigma_{\frak a_j}z\right)^s,$$
where the sum $\sum_*$ runs over all integers $a, d>0, 0\leq b<d, ad=n$
such that the matrix equation (3.8) has a solution $T,A,B,D$,
is equal to
$$\sum_{\underset{(w_j, N/w_j)|(A-D)}\to{A,D>0, AD=n}}
\{\delta_{ij}A^sD^{1-s}y^s+\varphi_{ij}(s)A^{1-s}D^sy^{1-s}\}.
\tag 3.9$$

    Next, we assume that $a, d>0$ and $b$ are integers with
$ad=n$ and $0\leq b<d$ such that the matrix equation (3.8)
has no solutions.   Then the cusp
$(a\frak a_j+b)/d$ is not $\Gamma_0(N)$-equivalent to $\frak a_j$.
Because every cusp of $\Gamma_0(N)$ is $\Gamma_0(N)$-equivalent to
a cusp given in (3.1), a cusp $\frak a_{j^\prime}$,
$j^\prime\neq j$, exists which is equivalent
to $(a\frak a_j+b)/d$ under $\Gamma_0(N)$.   That is
$$\frak a_{j^\prime}=\gamma_0\left({a\frak a_j+b\over d}\right)$$
for some element $\gamma_0$ in $\Gamma_0(N)$.
Let $\sigma_{\frak a_{j^\prime}}$ be the element in
$PSL_2(\Bbb R)$ such that $\sigma_{\frak a_{j^\prime}}\infty
=\frak a_{j^\prime}$ and $\sigma_{\frak a_{j^\prime}}^{-1}
\Gamma_{\frak a_{j^\prime}}\sigma_{\frak a_{j^\prime}}=\Gamma_\infty$,
and let
    $$T_0=\sigma_{\frak a_{j^\prime}}^{-1}\gamma_0{1\over\sqrt n}
\left(\matrix a&b\\0&d\endmatrix\right)\sigma_{\frak a_j}.$$
Then $T_0$ is an element of $SL_2(\Bbb R)$ satisfying
$T_0(\infty)=\infty$.    It follows that we have the decomposition
$$\Gamma_\infty\backslash\sigma_{\frak a_i}^{-1}
{1\over\sqrt n}\Gamma_0(N)\left(\matrix
a&b\\0&d\endmatrix\right)\sigma_{\frak a_j}
=\left[ \bigcup_{\gamma>0}
\bigcup_{\delta\,(\text{mod}\,\gamma)}
\left(\matrix *&*\\ \gamma &\delta
\endmatrix\right)\Gamma_\infty\right]T_0 $$
where the unions on $\gamma, \delta$ run over all numbers $\gamma>0$,
$0\leq \delta<\gamma$ such that
$\left(\smallmatrix *&*\\ \gamma&\delta
\endsmallmatrix\right)$ belong to $\sigma_{\frak a_i}^{-1}
\Gamma_0(N)\sigma_{\frak a_{j^\prime}}$.
Note that, since $f_i(\sigma_{\frak a_j} z)$ and
$\widetilde f_i(\sigma_{\frak a_j}z)$ are both periodic
function of $x$, $z=x+iy$, of period $1$, the function
$$\widetilde g_i(\sigma_{\frak a_j}z)
=f_i(\sigma_{\frak a_j} z)-\widetilde f_i(\sigma_{\frak a_j}z)$$
is also a periodic function of $x$ of period $1$.   Since
$T_0(\infty)=\infty$, we can write
$$T_0={1\over\sqrt n}\left(\matrix A^\prime &B^\prime
\\0&D^\prime\endmatrix\right)$$
for some real numbers $A^\prime,B^\prime$ and $D^\prime$.
 Then it follows the above decomposition, the above note,
  and the proof of Theorem 3.4 of Iwaniec [6]
that the constant coefficient in the Fourier
expansion at infinity of the function
$\widetilde g_i(\sigma_{\frak a_j} z)$ is equal to
$$\sum_{a,d,b}^*\varphi_{ij^\prime}(s)
\left({A^\prime\over D^\prime} y\right)^{1-s}\tag 3.10$$
where the sum is over all integers $a,d>0$ and $0\leq b<d$ with
$ad=n$ such that the matrix equation (3.8) has no solutions.
It follows from (3.7), (3.9) and (3.10) that
$$\aligned &\sum_{i=1}^{\nu(N)}
\int_{D_Y} f_i(z)E_i(z, \bar s)dz-o_Y(1)  \\
&=\sum_{i,j=1}^{\nu(N)}\sum_{\underset{(w_j, N/w_j)|(A-D)}
\to{A,D>0, AD=n}}\{{A^sD^{1-s}\delta_{ij}
Y^{s+\bar s-1}-\varphi_{ij}(s)A^{1-s}D^s
\varphi_{ij}(\bar s) Y^{-s-\bar s+1}\over s+\bar s-1}\\
&\hskip1.5truein
+{\delta_{ij}A^sD^{1-s}\varphi_{ij}(\bar s) Y^{s-\bar s}
-\varphi_{ij}(s)A^{1-s}D^s\delta_{ij}
 Y^{-s+\bar s}\over s-\bar s}\}\\
&-\sum_{i,j=1}^{\nu(N)}\sum_{a,d,b}\{{\left({A^\prime
\over D^\prime}\right)^{1-s}\varphi_{ij^\prime}(s)
\varphi_{ij}(\bar s) Y^{-s-\bar s+1}\over s+\bar s-1}+
{\delta_{ij}\left({A^\prime
\over D^\prime}\right)^{1-s}\varphi_{ij^\prime}(s)
 Y^{-s+\bar s}\over s-\bar s}\}\endaligned$$
where the sum over $a,d,b$ has the same meaning as in (3.10).
If $\Phi(s)$ is the constant term matrix of the Eisenstein
series for $\Gamma_0(N)$, then $\Phi(s)$ is a symmetric matrix
and satisfies the identity
$$\Phi(s)\Phi(1-s)=I$$
for all complex $s$, where $I$ is the identity matrix;
see Kubota [7].  Since $j\neq j^\prime$, we have
$$\sum_{i=1}^{\nu(N)}\varphi_{ij^\prime}(s)\varphi_{ij}(\bar s)=0$$
as $\text{Re}\, s\to 1/2$.   Hence, the above identity can
be written as
$$\aligned &\sum_{i=1}^{\nu(N)}
\int_{D_Y} f_i(z)E_i(z, \bar s)dz
=\sum_{i,j=1}^{\nu(N)}\sum_{\underset{(w_j, N/w_j)|(A-D)}
\to{A,D>0, AD=n}}A^sD^{1-s}\\
&\times \{{\delta_{ij}
Y^{s+\bar s-1}-\varphi_{ij}(s)
\varphi_{ij}(\bar s) Y^{-s-\bar s+1}\over s+\bar s-1}
+\delta_{ij}{\varphi_{ij}(\bar s) Y^{s-\bar s}
-\varphi_{ij}(s) Y^{-s+\bar s}\over s-\bar s}\}\\
&-\sum_{i,j=1}^{\nu(N)}\sum_{a,d,b}^*
{\delta_{ij}\left({A^\prime
\over D^\prime}\right)^{1-s}\varphi_{ij^\prime}(s)
 Y^{-s+\bar s}\over s-\bar s}+o_Y(1)\endaligned \tag 3.11$$
 as $\text{Re}\, s\to 1/2$, where the sum on $a,d,b$
 has the same meaning as in (3.10).

 This completes the proof of the theorem. \qed\enddemo

 \proclaim{Theorem 3.12}  Put
$$c(\infty)_Y=\delta_n\nu(N)\int_0^Y\int_0^1\sum_{0\neq b\in \Bbb Z}
k(z, z+\frac b{\sqrt n})dz -\int_{D_Y} H(z, z)dz.$$
Then
$$\aligned \frac {c(\infty)_Y}{\sqrt n} &=\delta_n\nu(N)
g(0)\ln{\sqrt n\over 2}-\sum_{j=1}^{\nu(N)}
\sum_{\underset{\frak a_j={u_j\over w_j}, (w_j, {N\over w_j})|(a-d)}
\to {ad=n, d>0, a\neq d}}g(\ln{a\over d})\ln Y\\
&-{\delta_n\nu(N)\over 2\pi}\int_{-\infty}^\infty
h(r)\frac{\Gamma^\prime}\Gamma(1+ir)dr +{1\over 4} h(0)
\{\delta_n\nu(N)\\
&+\sum_{j=1}^{\nu(N)}\sum_{\underset{\frak a_j={u_j\over w_j},
(w_j, {N\over w_j})|(a-d)}
\to {ad=n, d>0, a\neq d}}
\varphi_{jj}({1\over 2})\}\\
&+\sum_{i,j=1}^{\nu(N)}\sum_{\underset{
\frak a_j={u_j\over w_j},(w_j, {N\over w_j})|a-d}
\to {ad=n, d>0}}{1\over 4\pi}\int_{-\infty}^\infty h(r)
\left({a\over d}\right)^{ir}\varphi_{ij}^\prime({1\over 2}+ir)
\varphi_{ij}({1\over 2}-ir)dr\\
&-{1\over 8}h(0)\sum_{j=1}^{\nu(N)}\sum_{a,d,b}^*
\left({A^\prime\over D^\prime}\right)^{1/2}
\varphi_{jj^\prime}({1\over 2})+o_Y(1),\endaligned$$
where $o_Y(1)$ tends to zero as $Y\to\infty$, and
the sum on $a,d,b$ has the same meaning as in (3.10)
with $A^\prime, D^\prime$ being defined there.
\endproclaim

 \demo{Proof}  By the argument of [7], pp.102--106 we have
$$\aligned &\frac 1{\sqrt n}\int_0^Y \int_0^1
\sum_{0\neq b\in \Bbb Z}
k(z, z+\frac b{\sqrt n})dz \\
& =g(0)\ln(\sqrt n Y)-\frac 1{2\pi}\int_{-\infty}^\infty
h(r)\frac{\Gamma^\prime}\Gamma(1+ir)dr -g(0)\ln 2
+\frac 14 h(0)+o(1).\endaligned $$
Let
$$\varphi_{ij, m}(s)=\sum_c{1\over |c|^{2s}}\left(\sum_d
e(m{d\over c})\right)$$
where the sum runs over $c>0, d$ modulo $c$ with
$\left(\smallmatrix *&*\\ c&d\endsmallmatrix\right)\in
\sigma_{\frak a_i}^{-1}\Gamma_0(N)\sigma_{\frak a_j}$.
Then we have
$$\aligned E_i(\sigma_{\frak a_j}z, s)&=\delta_{ij}y^s+
\varphi_{ij}(s)y^{1-s}\\
&+{2\pi^s\sqrt y\over\Gamma(s)}\sum_{m\neq 0}
|m|^{s-{1\over 2}}K_{s-{1\over 2}}(2|m|\pi y)
\varphi_{ij, m}(s)e(mx)\endaligned$$
where
$$\varphi_{ij}(s)={\sqrt\pi\Gamma(s-{1\over 2})\over
\Gamma(s)}\varphi_{ij, 0}(s).$$
By Theorem 3.11, we have
$$\aligned &\sum_{i=1}^{\nu(N)}\int_{D_Y}
\{\sum_{ad=n, 0\leqslant b<d}E_i\left({az+b\over d},s\right)\}
E_i(z, \bar s)dz \\
&=\sum_{i,j=1}^{\nu(N)}\sum_{\underset{\frak a_j={u_j\over
w_j}, (w_j, N/w_j)|(a-d)}
\to{a,d>0, ad=n}}a^sd^{1-s} \{{\delta_{ij}
Y^{s+\bar s-1}-\varphi_{ij}(s)
\varphi_{ij}(\bar s) Y^{-s-\bar s+1}\over s+\bar s-1}\\
&+\delta_{ij}{\varphi_{ij}(\bar s) Y^{s-\bar s}
-\varphi_{ij}(s) Y^{-s+\bar s}\over s-\bar s}\}
-\sum_{i,j=1}^{\nu(N)}\sum_{a,d,b}^*
{\delta_{ij}\left({A^\prime
\over D^\prime}\right)^{1-s}\varphi_{ij^\prime}(s)
 Y^{-s+\bar s}\over s-\bar s}+o_Y(1)\endaligned $$
 as $\text{Re}\, s\to 1/2$, and hence we have
$$\aligned &\sum_{i=1}^{\nu(N)}\int_{D_Y}
\{\sum_{ad=n, 0\leqslant b<d}E_i\left({az+b\over d},s\right)\}
E_i(z, \bar s)dz \\
&=\sum_{j=1}^{\nu(N)}\sum_{\underset{\frak a_j=
{u_j\over w_j}, (w_j, {N\over w_j})|(a-d)}
\to{a,d>0, ad=n}}a^sd^{1-s} \{{Y^{s+\bar s-1}-\sum_{i=1}^{\nu(N)}
\varphi_{ij}(s)\varphi_{ij}(\bar s) Y^{-s-\bar s+1}\over s+\bar s-1}\\
&+{\varphi_{jj}(\bar s) Y^{s-\bar s}
-\varphi_{jj}(s) Y^{-s+\bar s}\over s-\bar s}\}
-\sum_{j=1}^{\nu(N)}\sum_{a,d,b}^*
{\left({A^\prime\over D^\prime}\right)^{1-s}\varphi_{jj^\prime}(s)
 Y^{-s+\bar s}\over s-\bar s}+o_Y(1)\endaligned $$
as $\text{Re}\, s\to 1/2$,
where the sum on $u/w$ runs over all cusps given in (3.1)
and the sum on $a,d,b$ has the same meaning as in (3.10)
with $A^\prime, D^\prime$ being defined there
and $o_Y(1)$ tends to zero as $Y\to\infty$.

By partial integration, we obtain
$$h(r)=\frac 1{r^4}\int_0^\infty g^{(4)}(\ln u) u^{ir-1}du $$
for nonzero $r$.   Then it follows that
$$\aligned &\lim_{S\to {1\over 2}^+}\int_{-\infty}^\infty
h(r)\{\sum_{i=1}^{\nu(N)}
\int_{D_Y} \left(\sum_{ad=n, 0\leqslant b<d}
E_i({az+b\over d}, S+ir)\right)E_i(z,S-ir)dz\}dr\\
&=\sqrt n\sum_{j=1}^{\nu(N)}\sum_{\underset{
\frak a_j={u_j\over w_j},(w_j, {N\over w_j})|a-d}
\to {ad=n, d>0}}\{4\pi g(\ln{a\over d})\ln Y
-\sum_{i=1}^{\nu(N)} \int_{-\infty}^\infty h(r)
\left({a\over d}\right)^{ir}\varphi_{ij}^\prime({1\over 2}+ir)\\
&\times \varphi_{ij}({1\over 2}-ir)dr +\int_{-\infty}^\infty h(r)
\left({a\over d}\right)^{ir}{\varphi_{jj}
({1\over 2}-ir) Y^{2ir}\over ir}dr\}\\
&-\sum_{j=1}^{\nu(N)}\sum_{a,d,b}^*
\int_{\infty}^\infty h(r)\varphi_{jj^\prime}({1\over 2}+ir)
\left({A^\prime\over D^\prime}\right)^{{1\over 2}-ir}
{Y^{-2ir}\over 2ir}dr+o_Y(1).\endaligned\tag 3.12$$
By the Riemann-Lebesgue theorem (cf. \S1.8 of [16]), we have
$$\aligned &\lim_{Y\to\infty}\int_{-\infty}^\infty h(r)
\left({a\over d}\right)^{ir}{\varphi_{jj}
({1\over 2}-ir) Y^{2ir}\over ir}dr\\
&=\lim_{Y\to\infty} \int_{-\infty}^\infty h(r)
\Re\left(\varphi_{jj}({1\over 2}-ir)\right)
{\sin\{r\ln(aY^2/d)\}\over r}dr=\pi h(0)
\varphi_{jj}({1\over 2}).\endaligned \tag 3.13$$
By (3.12) and (3.13), we have
$$\aligned &\frac 1{\sqrt n}\int_{D_Y}H(z, z)dz
= \sum_{j=1}^{\nu(N)}\sum_{\underset{
\frak a_j={u_j\over w_j},(w_j, {N\over w_j})|a-d}
\to {ad=n, d>0}}\{g(\ln{a\over d})\ln Y
+{1\over 4}h(0)\varphi_{jj}({1\over 2})\\
&-\sum_{i=1}^{\nu(N)}{1\over 4\pi}\int_{-\infty}^\infty h(r)
\left({a\over d}\right)^{ir}\varphi_{ij}^\prime({1\over 2}+ir)
\varphi_{ij}({1\over 2}-ir)dr\}\\
&+{1\over 8}h(0)\sum_{j=1}^{\nu(N)}\sum_{a,d,b}^*
\left({A^\prime\over D^\prime}\right)^{1/2}
\varphi_{jj^\prime}({1\over 2})+o_Y(1).\endaligned$$
 The stated identity then follows.  \qed\enddemo

   It follows from Theorem 3.9 and Theorem 3.12 that
$$\aligned &\lim_{Y\to\infty}\left( c(\infty)_Y
+\sum_{\{P\}, \Gamma_P=\{1_2\}}c(P)_Y\right)\\
&=\sqrt n\delta_n\nu(N)g(0)\ln{\sqrt n\over 2}+{\sqrt n\over 4}
h(0)\{\delta_n\nu(N)+\sum_{j=1}^{\nu(N)}
\sum_{\underset{\frak a_j={u_j\over w_j}, (w_j, {N\over w_j})|(a-d)}
\to {ad=n, d>0, a\neq d}}\varphi_{jj}({1\over 2})\}\\
&+{\sqrt n\over 2}\sum_{w|N, w>0}\sum_{ad=n, d>0, a\neq d}
\sum_b{\ln \{(a-d)^2wN[{C^2\over\ell^2},N]/C^2\}
\over |a-d|}g(\ln \frac ad)\\
&-{\delta_n\nu(N)\sqrt n\over 2\pi}\int_{-\infty}^\infty
h(r)\frac{\Gamma^\prime}\Gamma(1+ir)dr \\
&+\sum_{j=1}^{\nu(N)}\sum_{\underset{
\frak a_j={u_j\over w_j},(w_j, {N\over w_j})|a-d}
\to {ad=n, d>0}}
{\sqrt n\over 4\pi}\int_{-\infty}^\infty h(r)
\left({a\over d}\right)^{ir}\sum_{i=1}^{\nu(N)}
\varphi_{ij}^\prime({1\over 2}+ir)
\varphi_{ij}({1\over 2}-ir)dr\\
&+\sum_{\underset{\frak a_j={u_j\over w_j}, (w_j, {N\over w_j})|(a-d)}
\to {ad=n, d>0, a\neq d}}{1\over 2}
 |a-d|\int_1^\infty  k\left(\frac{(a-d)^2}nt\right)
\frac{\ln t}{\sqrt {t-1}}dt\\
 &-{1\over 8}h(0)\sum_{j=1}^{\nu(N)}\sum_{a,d,b}^*
\left({A^\prime\over D^\prime}\right)^{1/2}
\varphi_{jj^\prime}({1\over 2})\endaligned \tag 3.14$$
with $C=(a-d)w-bwN$, where $\ell=(C, bN)$ and the summation
on $b$ is taken over all numbers $b$ such that $C, bN/w\in \Bbb Z, N|C$
 and $0\leqslant b<|a-d|$.   Note that there are exactly $|a-d|$
number of such numbers $b$ by Lemma 3.5.

  Denote by $c(\infty)$ the right side of the identity (3.14).
We conclude that the trace formula (2.3) can be written as
$$\aligned & d(n)h(-\frac i2)+\sqrt n \sum_{j=1}^\infty h(\kappa_j)
\text{tr}_jT_n\\
&= c(I)+\sum_{\{R\}} c(R)+\sum_{\{P\},\,\Gamma_P\neq \{1_2\}}
c(P) +c(\infty) \endaligned \tag 3.15$$
for $\text{Re}\, s>1$, where the summations on the right side of
the identity are taken over the conjugacy classes.

   \proclaim{Theorem 3.13}   Let
$c(P)=\int_{\Gamma_P\backslash \Cal H} k(z, Pz)dz$
for hyperbolic elements $P\in\Gamma_0(N)$.  Then the series
$$\sum_{\{P\}} c(P)$$
represents an analytic function in the half-plane $\text{Re} s> 0$
except for having a possible pole at $s=1/2$ and for having
simple poles at $s=1, \frac 12\pm i\kappa_j$, $j=1,2,\cdots$.
 \endproclaim

  \demo{Proof}   We have
$$g^{(4)}(\log u)=A(s)u^{\frac 12-s}+O_s(u^{-\frac 12}),$$
where $A(s)$ is an analytic function of $s$
for Re$\, s>0$ and, for every complex number $s$ with
Re$\, s>0$, there exists a finite constant $B(s)$ depending
only on $s$ such that
$$|O_s(u^{-\frac 12})|\leqslant B(s)u^{-\frac 12}.$$
Moreover, for every fixed value of $u$, the term
$O_s(u^{-\frac 12})$ also represents an analytic
function of $s$ for Re$\, s>0$.  Since
$$h(r)=\frac 1{r^4}\int_0^\infty g^{(4)}(\ln u) u^{ir-1}du $$
for nonzero $r$, we have
$$ h(r)=\frac {A(s)}{r^4}\left(\frac 1{s-\frac 12-ir}+
\frac 1{s-\frac 12+ir}\right)+O_s(r^{-4}) $$
for Re$\, s>1$ and for nonzero $r$ with
$|\text{Im}\, r|<{1\over 2}-\epsilon$.
By analytic continuation, we obtain that
$$h(r)=\frac{2A(s)(s-\frac 12)}{r^4[(s-\frac 12)^2+
r^2]}+O_s(r^{-4}) \tag 3.16$$
for Re$\, s>0$ and for nonzero $r$ with
$|\text{Im} r|<{1\over 2}-\epsilon$.
It follows from results of [17] and (2.3)
that the left side of (3.15) is an analytic function of $s$
for Re$\, s>0$ except for having simple poles at
$s=1, \frac 12\pm i\kappa_j$,
$j=1,2,\cdots$.   Then the right
side of (3.15) can be interpreted as an analytic function of $s$ in
the same region by analytic continuation.

   Since $k(t)=(1+t/4)^{-s}$, by Lemma 3.1 we have that $c(R)$ is
analytic for Re$\, s>0$ except for a simple pole at $s=1/2$.
There are only a finite number of elliptic conjugacy classes
$\{R\}$.  The term $c(I)$ is a constant.

  Since $g(0)=2\sqrt\pi\Gamma(s-{1\over 2})\Gamma(s)^{-1}$,
$$h(0)=2\sqrt\pi 4^s{\Gamma(s-{1\over 2})\over\Gamma(s)}
\int_1^\infty (u+{1\over u}+2)^{{1\over 2}-s}{du\over u}$$
and
$$g(\ln{a\over d})=2\sqrt\pi 4^{s-{1\over 2}}{\Gamma(s-{1\over 2})
\over \Gamma(s)}\left({a\over d}+{d\over a}+2\right)^{{1\over 2}-s},$$
the sum of first three terms and the last two terms
on the right side of the identity (3.14)
is analytic for Re$\, s>0$ except for a pole at $s=1/2$.

 By Stirling's formula the identity
$${\Gamma^\prime(z)\over\Gamma(z)}=\ln z+O(1) \tag 3.17$$
holds uniformly when $|\arg z|\leq \pi-\delta$ for a small positive
number $\delta$.   It follows from (3.13) and (3.17) that the
fourth term on the right side of the identity (3.14)
is analytic for Re$\, s>0$ except for a possible pole at $s=1/2$.

   By Theorem 4.4.1 of Kubota \cite{7}, each Eisenstein series
$E_i(z,s)$ has a meromorphic continuation to the whole
$s$-plane, and the identity
$$\sum_{i=1}^{\nu(N)}\varphi_{ij}\left({1\over 2}+ir\right)
\varphi_{ij}\left({1\over 2}-ir\right)=1 $$
hold for all real $r$ and for $i=1,2,\cdots,\nu(N)$.
It follows that functions $\varphi_{ij}(s)$, $i,j=1,2,\cdots,\nu(N)$,
are analytic on the line $\text{Re}\,s=1/2$.
Let $Y$ be a fixed large positive number.  By the Maass-Selberg relation,
 we have
$$\aligned &\text{Re}\sum_{i=1}^{\nu(N)} \varphi_{ij}^\prime
\left({1\over 2} +ir\right) \varphi_{ij}\left({1\over 2}-ir\right)\\
&=2\ln Y +{\varphi_{jj}(1/2-ir)Y^{2ir}-\varphi_{jj}(1/2+ir)
Y^{-2ir}\over 2ir}\\
&-\int_{D_Y}E_j\left(z, {1\over 2}+ir\right)
E_j\left(z, {1\over 2}-ir\right)dz + o(1).
\endaligned \tag 3.18$$
Let $\frak a_i=u_i/w_i$ be a cusp given in (3.1),
and let $\eta=1/2+ir$.  By the argument following Lemma 3.6
of \cite {3}, we have
$$\aligned E_i(z, \eta)&=\delta_{\frak a_i\infty}y^\eta+
{\sqrt\pi\Gamma(\eta-{1\over 2})\over
\Gamma(\eta)}\varphi_{\frak a_i\infty, 0}(\eta)y^{1-\eta}\\
&+{2\pi^\eta\sqrt y\over\Gamma(\eta)}\sum_{m\neq 0}
|m|^{\eta-{1\over 2}}K_{\eta-{1\over 2}}(2|m|\pi y)
\varphi_{\frak a_i\infty, m}(\eta)e(mx)\endaligned \tag 3.19$$
where
$$\aligned &\varphi_{\frak a_i\infty, 0}(\eta)\\
&={\varphi(w_i)\over
\varphi((w_i,N/w_i))}\left({(w_i,N/w_i)\over w_iN}\right)^\eta
\prod_{p|N}{1\over 1-p^{-2\eta}}\prod_{p|{N\over w_i}}
(1-p^{1-2\eta}){\zeta{2\eta-1}\over\zeta(2\eta)}
\endaligned $$
and
$$\varphi_{\frak a_i\infty, m}(\eta)
=\left({(w_i,N/w_i)\over w_iN}\right)^\eta
\sum_{(c,N/w_i)=1}{1\over c^{2\eta}}
\sum_{\underset{cd\equiv u_i\,\text{mod} (w_i,N/w_i)}\to
{d\,\text{mod}(cw_i),\, (d,cw_i)=1}}
e\left(-{md\over cw_i}\right).$$
It follows from the functional identity of the Riemann
zeta-function that
$$\aligned {\Gamma(\eta-{1\over 2})\over \Gamma(\eta)}
\varphi_{\frak a_i\infty, 0}(\eta)
=&\pi^{2\eta-{3\over2}}{\Gamma(1-\eta)\over
\Gamma(\eta)}{\zeta(2-2\eta)\over\zeta(2\eta)}
{\varphi(w_i)\over \varphi((w_i,N/w_i))}\\
&\times \left({(w_i,N/w_i)\over w_iN}\right)^\eta
\prod_{p|N}{1\over 1-p^{-2\eta}}\prod_{p|{N\over w_i}}
(1-p^{1-2\eta}). \endaligned $$
By using Stirling's formula
$$|\Gamma(\sigma+it)|\sim \sqrt{2\pi}e^{-\pi |t|/2}|t|^{\sigma-1/2}$$
for any fixed real value of $\sigma$ as $t\to\infty$, we obtain that
$${\Gamma(\eta-{1\over 2})\over \Gamma(\eta)}
\varphi_{\frak a_i\infty, 0}(\eta)y^{1-\eta}
\ll \sqrt y \ln^2 (|r|+1)). \tag 3.20$$
By (1.1) and by partial integration, we find that
$$ {K_{\eta-1/2}(2|m|\pi y)\over\Gamma(\eta)}
={2^{\eta-1/2}\over (2|m|\pi y)^{\eta-1/2}\sqrt\pi}\int_0^\infty
{\cos (2|m|\pi yt)\over (1+t^2)^{\eta}}dt
\ll {1+|r|^3\over |m|^3 y^3}. \tag 3.21$$
Then it follows from (3.16), (3.18)-(3.21), Lemma 4.7 of \cite{3}
together with the proof of Theorem 1 in [8], and Cauchy's
inequality that the fifth term on the right side of the
identity (3.14) is analytic for Re$\, s>0$ except for a possible
pole at $s=1/2$.

  Also, it follows from (3.12) and Lemma 8.2 of Str\"ombergsson
[15] that  the fifth term on the right side of the
identity (3.14) is analytic for Re$\, s>0$ except for a possible
pole at $s=1/2$.

   Therefore, by (3.15) we have proved that  the series
$$\sum_{\{P\}} c(P)$$
represents an analytic function of $s$ in the half-plane
$\text{Re} s> 0$ except for having a possible pole
at $s=1/2$ and for having simple poles at $s=1, \frac 12\pm i\kappa_j$,
$j=1,2,\cdots$. \qed\enddemo

\heading
4.   Proof of the Main Theorem
\endheading

  Proofs for the following lemmas from Lemma 4.1 to Lemma 4.4
are minor modifications of corresponding lemmas
 in Conrey-Li [2].  For the convenience of readers,
we sketch the proofs here.

 \proclaim{Lemma 4.1}   Let $\lambda>0$ be an eigenvalue of
$\Delta$ for $\Gamma_0(N)$ with $N$ square free.  Assume that
$(n, N)=1$.  Put $\tau=1/2+i\kappa$ with $\kappa=\sqrt{\lambda-1/4}$.
Then we have
 $$ 4^\tau\sqrt{\pi n}\frac{\Gamma(i\kappa)}{\Gamma(\tau)}
 \text{tr}T_n=\lim_{s\to\tau}(s-\tau)
\sum_{\{P\}, \,\Gamma_P\neq \{1_2\}}c(P)$$ where the right side
is defined as in Theorem 3.13.
\endproclaim

\demo{Proof}  By (2.2), (3.13) and results of [17], we have
$$\lim_{s\to\tau}(s-\tau)\sum_{j=1, \kappa_j\neq\kappa}^\infty
h(\kappa_j)\text{tr}_jT_n=0.$$ By the proof of Theorem 3.13, we
have
$$\lim_{s\to\tau}(s-\tau)\left(d(n)h(-{i\over 2})
-c(I)-c(\infty)-\sum_{\{R\}}c(R)\right)=0.$$
  Then the stated identity follows from
(2.2), (3.12) and Theorem 3.13. \qed\enddemo

  \remark{Remark}  Let $h_d$ be the number of inequivalent classes of primitive
indefinite quadratic forms of discriminant $d$.  Siegel [14] proved that
$$\lim_{d\to\infty} \frac{\ln(h_d\ln\epsilon_d)}{\ln d}=\frac 12.\tag 4.1$$
 \endremark

  \proclaim{Lemma 4.2}   $P$ is a hyperbolic element of
$\Gamma^*$ with $\Gamma_P\neq\{1_2\}$ if and only if there
exists a primitive indefinite quadratic form $[a,b,c]$
of discriminant $d$ such that
 $$P=\frac 1{\sqrt n}\left(\matrix
{\frac {v-bNu(a, N)^{-1}}2}&{-cNu(a,N)^{-1}}\\
{aNu(a,N)^{-1}}&{\frac {v+bNu(a,N)^{-1}}2}\endmatrix\right)$$
with
$$v^2-\{dN^2(a,N)^{-2}\}u^2=4n.$$
   If   $\lambda_P$ is an eigenvalue of $P$, then
$$\lambda_P-{1\over\lambda_P}=\pm{Nu\over (a,N)}
{\sqrt d\over \sqrt n }.$$
Let
$$P_0=\left(\matrix {\frac {v_0-bu_0}2}&{-cu_0}\\
{au_0}&{\frac {v_0+bu_0}2}\endmatrix\right),$$
where the pair $(v_0, u_0)$ is the fundamental solution
of Pell's equation $v^2-du^2=4$.   Then $P$ is
$\Gamma_0(N)$-conjugate to a hyperbolic element
$P^\prime\in\Gamma^*$ with $\Gamma_{P^\prime}\neq \{1_2\}$
 if and only if $P_0$ is $\Gamma_0(N)$-conjugate to $P_0^\prime$,
where $P_0^\prime$ is associated with $P^\prime$ similarly as
$P_0$ is associated with $P$. \endproclaim

 \demo{Proof}   Let
$$P=\frac 1{\sqrt n}\left(\matrix A&B\\C&D\endmatrix\right)$$
be a hyperbolic element of $\Gamma^*$ such that $\Gamma_P\neq \{1_2\}$.
Then fixed points $r_1$, $r_2$ of $P$ are not rational numbers,
and hence $\Gamma_P$ is the subgroup of elements in $\Gamma_0(N)$
having $r_1$, $r_2$ as fixed points.
  Let $a=C/\mu, b=(D-A)/\mu$ and $c=-B/\mu$, where $\mu=(C, D-A, -B)$.
Then $[a, b, c]$ is a primitive quadratic form with
 $r_1$, $r_2$ being the roots of the equation $ar^2+br+c=0$.
By Sarnak [10], the subgroup of elements in $SL_2(\Bbb Z)$
 having $r_1$, $r_2$ as fixed points is generated by the
primitive hyperbolic element
$$P_0=\left(\matrix {\frac {v_0-bu_0}2}&{-cu_0}\\
{au_0}&{\frac {v_0+bu_0}2}\endmatrix\right)$$
where the pair $(v_0, u_0)$ is the fundamental solution of Pell's equation
$v^2-du^2=4$.

  Since $P$ and $P_0$ have the same fixed points, we have
$A=D-bC/a$ and $B=-cC/a$.  Since $P$ belongs to
$\Gamma^*$ and $AD-BC=n$, $C$ satisfies
$$\cases aD^2-bDC+cC^2=na\\ a|C,\,\, N|C.\endcases\tag 4.2$$
Let $\lambda_P$ be an eigenvalue of $P$.  Then it is a solution of
the equation $\lambda^2-{A+D\over\sqrt n}\lambda+1=0$.
By using $A=D-bC/a$ and $B=-cC/a$, we obtain that
$$\lambda_P-\frac 1{\lambda_P}=\pm
\frac {C\sqrt d}{a\sqrt n }  \tag 4.3$$
and
$$\lambda_P+\frac 1{\lambda_P}=\frac 1{\sqrt n}(2D-\frac ba C). \tag 4.4$$

  Conversely, let a pair $(C, D)$ be a solution of the equation (4.2).
Define $A=D-bC/a$ and $B=-cC/a$.  Then the matrix
$$P=\frac 1{\sqrt n}\left(\matrix A&B\\C&D\endmatrix\right)$$
is a hyperbolic element of $\Gamma^*$
with $\Gamma_P\neq\{1_2\}$ having the same fixed points as $P_0$,
and eigenvalues of $P$ satisfies (4.3) and (4.4).

  Next, let $v=2D-b\frac Ca$ and $u_1=\frac Ca$.
Then the equation (4.2) becomes $v^2-du_1^2=4n $ with $N|au_1$.
Since $N|au_1$, this equation can be written as
$$v^2-{dN^2\over (a,N)^2}u^2=4n. \tag 4.5$$
Moreover, we have
$$\frac 1{\sqrt n}\left(\matrix A&B\\C&D\endmatrix\right)
=\frac 1{\sqrt n}\left(\matrix
{\frac {v-bNu(a, N)^{-1}}2}&{-cNu(a,N)^{-1}}\\
{aNu(a,N)^{-1}}&{\frac {v+bNu(a,N)^{-1}}2}
\endmatrix\right). \tag 4.6$$

   Let $P_0^\prime$ be the primitive hyperbolic element of
$SL_2(\Bbb Z)$ corresponding to $[a^\prime, b^\prime, c^\prime]$.
Two forms $[a, b, c]$ and $[a^\prime, b^\prime, c^\prime]$
of the same discriminant are equivalent in $\Gamma_0(N)$
if and only if an element $\gamma\in\Gamma_0(N)$ exists
such that $\gamma^{-1}P_0\gamma=P_0^\prime$. \qed\enddemo

  \proclaim{Lemma 4.3}  Let $[a,b,c]$ be a primitive
indefinite quadratic form of discriminant $d$, and let
 $$P=\frac 1{\sqrt n}\left(\matrix
{\frac {v-bNu(a, N)^{-1}}2}&{-cNu(a,N)^{-1}}\\
{aNu(a,N)^{-1}}&{\frac {v+bNu(a,N)^{-1}}2}\endmatrix\right)$$
with $v^2-\{dN^2(a,N)^{-2}\}u^2=4n$
be a hyperbolic element of $\Gamma^*$
with $\Gamma_P\neq \{1_2\}$.  Let $(v_1, u_1)$
with $v_1, u_1>0$ be the
fundamental solution of Pell's equation
$v^2-d_1u^2=4$, where $d_1=dN^2/(a, N)^2$.
Then $\Gamma_P$ is generated by the hyperbolic element
$$P_1=\left(\matrix {{1\over 2}\left(v_1-b{Nu_1\over (a, N)}
\right)}&{-c{Nu_1\over (a, N)}}\\
{a{Nu_1\over (a, N)}}&{{1\over 2}\left(v_1+b{Nu_1\over
(a, N)} \right)}\endmatrix\right)$$
of $\Gamma_0(N)$.\endproclaim

  \demo{Proof}   Let $a, b$ and $c$ be given as in the proof
of Lemma 4.2.  Since $\Gamma_P$ is cyclic, a solution $v_1, u_1^\prime>0$
of Pell's equation $v^2-du^2=4$ exists such that
 $\Gamma_P$ is generated by
$$P_1=\left(\matrix {\frac {v_1-bu_1^\prime}2}&{-cu_1^\prime}\\
{au_1^\prime}&{\frac {v_1+bu_1^\prime}2}\endmatrix\right),$$
and hence $P_1$ is the smallest positive integer
power of $P_0$ among all powers of $P_0$ belonging to
$\Gamma_0(N)$.   Since $N|au_1^\prime$,
we have ${N\over (a, N)}|u_1^\prime$.  Write
$$u_1^\prime={Nu_1\over (a, N)}.$$
Then the pair $(v_1, u_1)$ with $v_1, u_1>0$ is
 the fundamental solution of the Pell equation
$v^2-d_1u^2=4$, where $d_1=dN^2(a, N)^{-2}$.
The stated result then follows. \qed\enddemo

 Two quadratic forms $[a, b,c]$ and $[a^\prime, b^\prime,
c^\prime]$ are equivalent under $\Gamma_0(N)$ if an element
$\gamma\in \Gamma_0(N)$ exists such that
$$\left(\matrix {a^\prime}&{b^\prime/2}\\
{b^\prime/2}&{c^\prime}\endmatrix\right)=\gamma^t
\left(\matrix a&{b/2}\\{b/2}&c\endmatrix\right)\gamma.$$

  \proclaim{Lemma 4.4}   Let $[a_j,b_j,c_j]$, $j=1,2,\cdots,H_d$,
be a set of representatives for classes of primitive indefinite
quadratic forms of discriminant $d$, which are not equivalent
under $\Gamma_0(N)$.  Then we have
$$\aligned &\sum_{\{P\}:\,\,P\in\Gamma^*,\,
\Gamma_P\neq \{1_2\}}c(P)\\
&= 4\sqrt{\pi n}{\Gamma(s-{1\over 2})\over N\Gamma(s)}
\sum_{d\in \Omega}\sum_{j=1}^{H_d}\sum_u
{(a_j,N)\ln \epsilon_{d_1}\over u\sqrt d}
\left(1+{d(Nu)^2\over 4n(a_j,N)^2}\right)^{{1\over 2}-s}
\endaligned$$
for $\Re s>1$, where $d_1=dN^2(a_j, N)^{-2}$ and
 the sum on $u$ is over all the
positive integers $u$ such that $4n+dN^2(a_j,N)^{-2}u^2$ is
the square of an integer.  \endproclaim

 \demo{Proof}  It follows from (2.1) and Lemma 3.2 that
$$\sum_{\{P\}:\,\,P\in\Gamma^*,\, \Gamma_P\neq \{1_2\}}c(P)
=2\sqrt\pi{\Gamma(s-{1\over 2})\over\Gamma(s)}\sum_{\{P\}}
{\ln NP_1\over \lambda_P-1/\lambda_P}\left(1+{(\lambda_P
-1/\lambda_P)^2\over 4}\right)^{{1\over 2}-s}$$
for $\Re s>1$, where $\lambda_P>1$ is an eigenvalue of $P$
and $P_1$ is given in Lemma 4.3.  Let $P$ be associated with
a primitive indefinite quadratic form $[a_j, b_j,c_j]$
as in Lemma 4.2.  Then by Lemma 4.2, we have
$$\lambda_P-{1\over\lambda_P}={Nu\over (a_j,N)}
{\sqrt d\over \sqrt n }.$$
By Lemma 4.3, we have
$$\sqrt{NP_1}={v_1+\sqrt{d_1}u_1\over 2}=\epsilon_{d_1}.$$

  If $P^\prime$ is a hyperbolic of $\Gamma^*$ with
$\Gamma_{P^\prime}\neq\{1_2\}$, and is associated with
a primitive indefinite quadratic form
$[a_{j^\prime}, b_{j^\prime},c_{j^\prime}]$
as in Lemma 4.2, then $P$ and $P^\prime$ are
$\Gamma_0(N)$-conjugate if and only if
$[a_j, b_j,c_j]$ and $[a_{j^\prime}, b_{j^\prime},c_{j^\prime}]$
are $\Gamma_0(N)$-conjugate.
  Let $T$ be a hyperbolic of $\Gamma^*$ with
$\Gamma_T\neq\{1_2\}$ associated with
a primitive indefinite quadratic form $[a,b,c]$
as in Lemma 4.2.  If the discriminant of $[a,b,c]$ is
not equal to the discriminant of $[a_j, b_j,c_j]$
which is associated with $P$, then $P$ and $T$ are
not $\Gamma_0(N)$-conjugate.  The stated identity
then follows.  \qed\enddemo

  The proof for the following lemma is a minor modification
of a corresponding lemma in Conrey-Li [2], and is given in
Cardon-Li [1].

  \proclaim{Lemma 4.5}  Let $k$ be a divisor of $N$.  Then
the number of indefinite primitive quadratic forms
$[a,b,c]$ with $(a,N)=k$ of discriminant $d$, which are not
equivalent under $\Gamma_0(N)$, is equal to
$$h_{d_1}\prod_{p^{2l}|(d,k)}p^l
\cdot\prod_{p|k}\left(1+\left({d\over p}\right)\right)$$
where $d_1=dN^2/k^2$, where the first product runs over all
distinct primes $p$ with $p^{2l}$ being the greatest even
$p$-power factor of $(d,k)$ and the second product is taken over all
distinct primes $p$ dividing $k$.  \endproclaim

    \proclaim{Theorem 4.6}   Let $N$ be a non-square free positive integer
with $(n, N)=1$.   Then we have
$$\aligned &\sum_{\{P\}:\,\,P\in\Gamma^*,\,
\Gamma_P\neq \{1_2\}}c(P)
= 4\sqrt{\pi n}{\Gamma(s-{1\over 2})
\over \Gamma(s)}\\
 &\times \sum_{k|N}\sum_{m|{N\over k}}\sum_{n|k}{\mu(n)\over n}
\sum_{d\in \Omega}\sum_u \prod_{p^{2l}|(d, {N\over k})}p^l
\left({d\over mn}\right)
{h_d\ln \epsilon_d\over u\sqrt d}
\left(1+{d k^2u^2\over 4n}\right)^{{1\over 2}-s}\endaligned $$
for $\Re s>1$, where the sum on $u$ is over
all the positive integers $u$ such that
$\sqrt{4n+d k^2 u^2}\in\Bbb Z$.\endproclaim

  \demo{Proof}  Let $[a_j,b_j,c_j]$, $j=1,2,\cdots,H_d$,
be a set of representatives for classes of primitive indefinite
quadratic forms of discriminant $d$, which are not equivalent
under $\Gamma_0(N)$.  By Lemma 4.4, we have
$$\sum_{\{P\}:\,\,P\in\Gamma^*,\,
\Gamma_P\neq \{1_2\}}c(P)
= 4\sqrt{\pi n}{\Gamma(s-{1\over 2})\over \Gamma(s)}
\sum_{d\in \Omega}\sum_{j=1}^{H_d}\sum_u
{\ln \epsilon_{d_1}\over u\sqrt {d_1}}
\left(1+{d_1 u^2\over 4n}\right)^{{1\over 2}-s}$$
for $\Re s>1$, where $d_1=dN^2(a_j, N)^{-2}$ and
the sum on $u$ is over all the
positive integers $u$ such that $\sqrt{4n+d_1 u^2}
\in\Bbb Z$.   Denote $k=N/(a_j, N)$.
By  Lemma 4.5, we can write the above identity as
$$\aligned &\sum_{\{P\}:\,\,P\in\Gamma^*,\,
\Gamma_P\neq \{1_2\}}c(P)
= 4\sqrt{\pi n}{\Gamma(s-{1\over 2})\over \Gamma(s)}\\
&\times\sum_{k|N}\sum_{d\in \Omega}\sum_u
\prod_{p^{2l}|(d, {N\over k})}p^l
\prod_{p|{N\over k}}\left(1+\left({d\over p}\right)
\right){h_{d_1}\ln \epsilon_{d_1}\over u\sqrt {d_1}}
\left(1+{d_1u^2\over 4n}\right)^{{1\over 2}-s}
\endaligned \tag 4.6$$
for $\Re s>1$, where $d_1=dk^2$.   By Dirichlet's
class number formula
$$h_{d_1}\ln\epsilon_{d_1}=\sqrt{d_1}L(1, \chi_{d_1})$$
and the identity
$$L(1, \chi_{d_1})=L(1, \chi_d)\prod_{p|k}
\left(1-\left({d\over p}\right)p^{-1}\right),$$
we can write (4.6) as
$$\aligned &\sum_{\{P\}:\,\,P\in\Gamma^*,\,
\Gamma_P\neq \{1_2\}}c(P)= 4\sqrt{\pi n}{\Gamma(s-{1\over 2})
\over \Gamma(s)}\sum_{k|N}\sum_{d\in \Omega}\sum_u
\prod_{p^{2l}|(d, {N\over k})}p^l\\
&\times \prod_{p|{N\over k}}\left(1+\left({d\over p}\right)
\right)\prod_{p|k}\left(1-\left({d\over p}\right)p^{-1}\right)
{h_d\ln \epsilon_d\over u\sqrt d}
\left(1+{d k^2u^2\over 4n}\right)^{{1\over 2}-s}\endaligned $$
for $\Re s>1$, where the sum on $u$ is over
all the positive integers $u$ such that
$\sqrt{4n+d k^2 u^2}\in\Bbb Z$.  Since
$$\prod_{p|{N\over k}}\left(1+\left({d\over p}\right)
\right)\prod_{p|k}\left(1-\left({d\over p}\right)p^{-1}\right)
=\sum_{m|{N\over k}}\sum_{n|k} \left({d\over mn}\right){\mu(n)\over n},$$
we have
$$\aligned &\sum_{\{P\}:\,\,P\in\Gamma^*,\,
\Gamma_P\neq \{1_2\}}c(P)= 4\sqrt{\pi n}{\Gamma(s-{1\over 2})
\over \Gamma(s)}\\
&\times \sum_{k|N}\sum_{m|{N\over k}}\sum_{n|k}{\mu(n)\over n}
\sum_{d\in \Omega}\sum_u \prod_{p^{2l}|(d, {N\over k})}p^l
\left({d\over mn}\right){h_d\ln \epsilon_d\over u\sqrt d}
\left(1+{d k^2u^2\over 4n}\right)^{{1\over 2}-s}.\endaligned
\tag 4.7$$

  Next, we want to show that
$$\sum_{d\in \Omega}\sum_u \left({d\over mn}\right)
\prod_{p^{2l}|(d, {N\over k})}p^l{h_d\ln \epsilon_d\over u\sqrt d}
\left(1+{d k^2u^2\over 4n}\right)^{{1\over 2}-s}$$
is absolutely convergent for $\sigma=\Re s>1$.  Since
$$\aligned & \left |\sum_{d\in \Omega}\sum_u
\prod_{p^{2l}|(d, {N\over k})}p^l
\left({d\over mn}\right){h_d\ln \epsilon_d\over u\sqrt d}
\left(1+{d k^2u^2\over 4n}\right)^{{1\over 2}-s}\right|\\
&\leqslant \sqrt N\sum_{d\in \Omega}\sum_u
{h_d\ln \epsilon_d\over u\sqrt d}
\left(1+{d k^2u^2\over 4n}\right)^{{1\over 2}-\sigma}.
\endaligned \tag 4.8$$
It is proved in Li [9] that the right side of (4.8)
is convergent for $\sigma>1$, and hence, the right side of
the stated identity is absolutely convergent for
$\Re s>1$.

This completes the proof of the theorem. \qed\enddemo

\vskip0.15truein

 \demo{Proof of Theorem 1}  It is proved at the end of Li [9] that
$$\sum_{d\in\Omega, u} \left|\frac{h_d\ln \epsilon_d} {\sqrt d u}
(1+\frac {du^2}{4n})^{\frac 12-\sigma}-(4n)^{\sigma-{1\over 2}}
\frac{h_d\ln \epsilon_d} {(du^2)^\sigma}\right|\ll \sum_{d\in\Omega, u}
\frac{(du^2)^{\frac 12+\epsilon-1-\sigma}}{u^{1+2\epsilon}}<\infty$$
for $\sigma>0$.  Then it follows from Theorem 4.6 that
$$\aligned & \lim_{s\to 1/2+i\kappa}(s-{1\over 2}-i\kappa)\sum_{\{P\}, \,
\Gamma_P\neq \{1_2\}}c(P)=\frac{4\sqrt{n\pi}\Gamma(i\kappa)(4n)^{i\kappa}}
{\Gamma(1/2+i\kappa)}\\
&\times \sum_{k|N}\sum_{m|{N\over k}}\sum_{n|k}{\mu(n)\over n}k^{-2i\kappa}
\lim_{s\to 1/2+i\kappa}(s-{1\over 2}-i\kappa)\sum_{d\in \Omega}\sum_u
\prod_{p^{2l}|(d, {N\over k})}p^l
\left({d\over mn}\right){h_d\ln \epsilon_d\over (du^2)^s}.
\endaligned \tag 4.9 $$
Theorem 3.13 shows that the function on the right side of (4.9)
represents an analytic function in the half-plane $\Re s> 0$
except for a possible pole at $s=1/2$ and for possible simple
poles at $s=1, \frac 12\pm i\kappa_j$, $j=1,2,\cdots$.  The
stated identity then follows from Lemma 4.1.

  This completes the proof of the theorem.  \enddemo

\enddocument